\documentclass[11pt]{amsart}
\usepackage[latin1]{inputenc}
\usepackage[T1]{fontenc}
\usepackage[english,francais]{babel}
\usepackage{ae, aecompl}
\usepackage{amsmath,amsfonts,amssymb}
\usepackage{amsthm}
\usepackage{amscd}
\usepackage{mathrsfs,array,graphicx}
\usepackage{stmaryrd}
\usepackage{wasysym}
\usepackage[all,ps]{xy}
\usepackage{calc}
\usepackage{enumerate}
\usepackage{ulem}

\usepackage{amscd}
\input xy
\xyoption{all}

\usepackage{fancyhdr}
\addto\captionsfrench{}

\newcommand{\ugu}{\stackrel{\mathrm{d\acute{e}f}}{=}}

\newcommand{\Q}{\mathbb{Q}_p}

\newcommand{\OF}{\mathcal{O}_F}
\newcommand{\Oe}{\mathcal{O}_E}
\newcommand{\he}{\mathcal{H}}
\newcommand{\Ind}{\mathrm{Ind}}
\newcommand{\ind}{\mathrm{c}\textrm{-}\mathrm{Ind}_{KZ}^G}
\newcommand{\into}{\hookrightarrow}
\newcommand{\onto}{\twoheadrightarrow}
\newcommand{\GL}{\mathrm{GL}_2}

\newcommand{\Sym}{\mathrm{Sym}}
\newcommand{\uni}{\varpi_F}

\newtheoremstyle{Proposition}{11pt}{11pt}{\itshape}{}{\bfseries}{.}{.5em }{}
\theoremstyle{Proposition}

\newcounter{Proposition}
\newtheorem{theo}[Proposition]{Théorème}
\newtheorem{prop}[Proposition]{Proposition}

\newtheorem{defin}[Proposition]{Définition}

\newtheorem{conj}[Proposition]{Conjecture}
\newtheorem{rem}[Proposition]{Remarque}
\newtheorem{cor}[Proposition]{Corollaire}
\newtheorem{lemma}[Proposition]{Lemme}

\numberwithin{equation}{section}
\numberwithin{Proposition}{section}
\author{Marco De Ieso}
\address{Bâtiment 430, Départment de Mathématiques d'Orsay, Université Paris-Sud, 91405 Orsay Cedex, France}
\email{Marco.DeIeso@math.u-psud.fr}

\title{EXISTENCE DE NORMES INVARIANTES POUR $\mathrm{GL}_2$}

\begin{document}
\date{}

\pagestyle{fancy}


\begin{otherlanguage}{english}
\begin{abstract}
In \cite{bs} Breuil et Schneider formulated  a conjecture on the equivalence of the existence of invariant norms on certain locally algebraic representations of $\mathrm{GL}_d(F)$ and the existence of certain de Rham representations of $\mathrm{Gal}(\overline{\mathbb{Q}}_p/F)$, where $F$ is a finite extension of $\Q$. In this paper we prove that in the case $d = 2$ and under some conditions, the existence of certains admissible filtrations implies the existence of invariant norms.

\vskip.5cm
  \par\noindent \textsc{R{\'e}sum{\'e}.} Dans \cite{bs} Breuil et Schneider formulent une conjecture sur l'equivalence entre l'existence de normes invariantes sur certaines représentations localement algébriques de $\mathrm{GL}_d(F)$ et l'existence de certaines représentations de de Rham de $\mathrm{Gal}(\overline{\mathbb{Q}}_p/F)$, où $F$ est une extension finie de $\Q$. Dans ce papier nous montrons que dans le cas  $d=2$ et sous  certaines  conditions, l'existence de certaines filtrations admissibles implique l'existence de normes invariantes. 
\end{abstract}
\end{otherlanguage}

\maketitle

\tableofcontents
\addtocontents{toc}{\protect\setcounter{tocdepth}{1}}

\section{Introduction, notations et énoncé des résultats} 
\subsection{Introduction}

Soit $p$ un nombre premier et $F$ une extension finie de $\Q$. Cet article s'inscrit dans le cadre du programme de Langlands local $p$-adique, qui a pour objet de relier certaines représentations $p$-adiques continues de dimension $d$ de $\mathrm{Gal}(\overline{\mathbb{Q}}_p/F)$ avec certaines représentations de $\mathrm{GL}_d(F)$.
\medbreak
Si $F = \Q$ et $d = 2$ alors tout est essentiellement bien compris: on dispose à présent d'une correspondance $V \mapsto \Pi(V)$ (\cite{colmez2}, \cite{pask})  associant à une $E$-représentation $V$ de $\mathrm{Gal}(\overline{\mathbb{Q}}_p/\Q)$, de dimension $2$, une représentation unitaire admissible de $\GL(\Q)$. Cette correspondance est compatible avec la correspondance de Langlands locale classique et avec la cohomologie étale complétée (\cite{eme2}).  
\medbreak
Les autres cas s'annoncent beaucoup plus délicats. En particulier, Breuil et Schneider ont formulé dans \cite{bs} une conjecture (géneralisant une conjecture antérieure de Schneider et Teitelbaum \cite{st}) qui laisse entrevoir un lien profond entre la catégorie des représentations continues de $\mathrm{Gal}(\overline{\mathbb{Q}}_p/\Q)$ de dimension $d$ et qui sont de de Rham, et certaines représentations localement algébriques de $\mathrm{GL}_d(F)$. L'idée à l'origine de cette conjecture est la suivante. D'après la théorie de Colmez et Fontaine (\cite{cf}) on sait qu'une représentation de de Rham peut être décrite par un espace vectoriel muni d'une action du groupe de Weil-Deligne de $F$ et d'une filtration, les deux étant reliés par une relation appelée faible admissibilité. À cet objet on peut associer une représentation lisse $\pi$ de $\mathrm{GL}_d(F)$ par la correspondance de Langlands modifiée. D'autre part, les poids de Hodge-Tate de la filtration permettent de construire une représentation algébrique irréductible de $\mathrm{GL}_d(F)$ que l'on désigne par $\rho$. La conjecture de Breuil et Schndeider dit alors essentiellement que l'existence d'une filtration faiblement admissible devrait être équivalente à l'existence d'une norme sur la représentation localement algébrique $\rho \otimes \pi$. Mentionnons que des résultats partiels ont été obtenus par Hu (\cite{hu}) et Sorensen (\cite{sore}).   
\medbreak
Dans cet article nous explorons le cas particulier où $F$ est une extension finie quelconque de $\Q$, $d=2$ et la représentation galoisienne de départ est cristalline. Avec, de plus, des hypothèses techniques supplémentaires sur les poids de la filtration nous donnons une réponse positive à cette conjecture. Les techniques que nous employons pour démontrer ce résultat sont classiques (\cite{breuilab}).






\subsection{Notations} Soit $p$ un nombre premier. On fixe une clôture algébrique $\overline{\mathbb{Q}}_p$ de $\Q$ et une extension finie $F$ de $\Q$ contenue dans $\overline{\mathbb{Q}}_p$. On note $\mathcal{O}_F$ l'anneau des entiers de $F$, $\mathfrak{p}_F$ son idéal maximal et $k_F = \mathcal{O}_F/\mathfrak{p}_F$ son corps résiduel. On note $f = [k_F:\mathbb{F}_p]$, $q = p^f$ et $e$ l'indice de ramification de $F$ sur $\Q$, de sorte que $[F:\Q] = ef$ et $k_F \simeq \mathbb{F}_q$. On note $F_0 = Frac(W(\mathbb{F}_q))$ le sous-corps non-ramifié maximal de $F$ et $\varphi_0$ le Frobenius sur $F_0$. On note $\mathrm{Gal}(\overline{\mathbb{Q}}_p/F)$ le groupe de Galois de $F$ et $W(\overline{\mathbb{Q}}_p/F)$ son groupe de Weil. La théorie du corps de classes local fournit une application $\mathrm{rec}\colon W(\overline{\mathbb{Q}}_p/F)^{ab} \to F^{\times}$ que l'on normalise en décidant d'envoyer les Frobenius arithmétique sur les inverses des uniformisantes. 
\medbreak 

On fixe une uniformisante $\uni \in \mathfrak{p}_F$ et on note $val_F$ la valuation $p$-adique sur $\overline{\mathbb{Q}}_p$ normalisée par $val_F(p) = [F:\Q]$. Si $x \in \overline{\mathbb{Q}}_p$ on pose $|x| = p^{- val_F(x)}$. Si $\lambda \in k_F$ alors $[\lambda]$ désigne le représentant de Teichmüller de $\lambda$ dans $\OF$. Si $\mu \in E^{\times}$ on note $\mathrm{nr}(\mu)\colon F^{\times} \to E^{\times}$ le caractère défini par $x \mapsto \mu^{val_F(x)}$.
\medbreak

On désigne par $G$ le groupe $\GL(F)$, par $K$ le groupe $\GL(\OF)$ qui est, à conjugaison près, l'unique sous-groupe compact maximal de $G$, par $I$ le sous-groupe d'Iwahori de $K$ et par $I(1)$ son pro-$p$-Iwahori. L'application de reduction modulo $\mathfrak{p}_F$ induit un homomorphisme surjectif:
\[
red: K \longrightarrow \GL(k_F).
\]
On rappelle que $I$ est l'ensemble des éléments de $K$ dont la réduction modulo $\mathfrak{p}_F$ est une matrice triangulaire supérieure et que $I(1)$ est le sous-groupe des éléments de $I$ dont la réduction modulo $\mathfrak{p}_F$ est une matrice unipotente. On désigne par $Z \simeq F^{\times}$ le centre de $G$ et par $P$ le sous-groupe de Borel formé des matrices triangulaires supérieures de $G$. On note:
\[
\alpha =  \begin{bmatrix} {1} & {0} \cr {0} & {\uni} \end{bmatrix}, \quad w = \begin{bmatrix} {0} & {1} \cr {1} & {0} \end{bmatrix}, \quad \beta = \alpha w = \begin{bmatrix} {0} & {1} \cr {\uni} & {0} \end{bmatrix}
\]
et, si $\lambda \in \OF$,
\[
w_{\lambda} = \begin{bmatrix} {0} & {1} \cr {1} & {-\lambda} \end{bmatrix}.
\]
\medbreak
On désignera toujours par $E$ une extension finie de $\Q$ qui vérifie
\[
|S| = [F:\Q],
\]
où $S \ugu \mathrm{Hom}_{alg}(F,E)$.

Si $\vec{n} = (n_{\sigma})_{\sigma \in S}, \vec{m} = (m_{\sigma})_{\sigma \in S}$ sont des $|S|$-uplets d'entiers positifs ou nuls posons: 
\begin{itemize}
\item[(i)] $\vec{n}! = \prod_{\sigma \in S}n_{\sigma}!$;
\item[(ii)] $|\vec{n}| = \sum_{\sigma \in S}n_{\sigma}$;
\item[(iii)] $\vec{n}-{\vec{m}} = (n_{\sigma}-m_{\sigma})_{\sigma \in S}$; 
\item[(iv)] $\vec{n}\leqslant{\vec{m}}$ si $n_{\sigma}\leq m_{\sigma}$ pour tout $\sigma \in S$;
\item[(v)] $\binom{\vec{n}}{\vec{m}} = \frac{\vec{n}!}{\vec{m}!(\vec{n}-\vec{m})!}$.
\end{itemize} 
Si $\vec{n}= (n_{\sigma})_{\sigma \in S} \in \mathbb{Z}_{\geq 0}^{|S|}$ et $z \in \OF$ on pose $z^{\vec{n}} = \prod_{\sigma \in S}\sigma(z)^{n_{\sigma}}$.

\subsection{Énoncé des résultats}

On fixe
\smallbreak
\begin{itemize}
\item[$\bullet$] $(\alpha,\beta) \in E^{\times} \times E^{\times}$.
\item[$\bullet$] un $|S|$-uplet d'entiers positifs $\vec{d}$. 
\end{itemize}
\medbreak
Notons:
\[
S^+ = \{\sigma \in S, d_{\sigma} \neq 0 \}\subseteq S 
\]  
et, pour $l$ dans $\{0,\ldots, f-1 \}$ posons:
\[
J_l = \Big\{\sigma\in S^+, \sigma([\zeta]) = \iota([\zeta])^{p^l} \ \forall \zeta \in k_F  \Big\}.
\]
Par exemple, si $F$ est non ramifiée alors $|J_l| = 1 $ pour tout $l$. 
\medbreak
Si $i \in \mathbb{Z}$ on désigne par $\overline{i}$ l'unique élément de $\{1,\ldots,f \}$ tel que $i \equiv \overline{i}\, \mathrm{mod}\, f$ et posons pour $\sigma \in J_l$:
\[
v_{\sigma} = \inf \Big\{1\leq i \leq f, J_{\overline{l+i}} \neq \emptyset  \Big\}.
\]

Notons $(\mathrm{Sym}^{d_{\sigma}}E^2)^{\sigma}$, pour $\sigma \in S$ et $d_{\sigma} \in \mathbb{Z}_{\geq 0}$, la représentation algébrique irréductible de $\mathrm{GL}_2 \otimes_{F,\sigma} E$ dont le plus haut poids est $\chi_{\sigma}\colon \mathrm{diag}(x_1,x_2) \mapsto \sigma(x_2)^{d_{\sigma}}$ vis-à-vis du sous-groupe des matrices triangulaires supérieures et posons:
\[
\rho_{\vec{d}} = \otimes_{\sigma \in S} (\mathrm{Sym}^{d_{\sigma}}E^2)^{\sigma}
\]
Notons
\[
\pi = \Ind_P^G (\mathrm{nr}(\alpha^{-1})\otimes \mathrm{nr}(p\beta^{-1})).
\]
l'induite parabolique lisse et non ramifiée usuelle. 
\medbreak
La conjecture de Breuil et Schneider dans ce cas peut être reformulée comme suit (voir la Section \ref{conseguenz} pour plus de détails).

\begin{conj}\label{congetturaintro}
Les deux conditions suivantes sont équivalentes:
\medbreak
\begin{itemize}
\item[(i)] La représentation $\rho_{\vec{d}}\otimes \pi$ admet une norme $G$-invariante, i.e. une norme $p$-adique telle que $\|gv\| = \|v\|$ pour tout $g\in G$ et $v \in \rho_{\vec{d}}\otimes \pi$.
\item[(ii)]  les inégalités suivantes sont vérifiées:
\medbreak
\begin{itemize}
\item[(1)] $val_F(\alpha^{-1})+val_F(p\beta^{-1})+\sum_{\sigma \in S}d_{\sigma} = 0$; 
\item[(2)] $val_F(p\beta^{-1})+\sum_{\sigma \in S}d_{\sigma} \geq 0$.
\end{itemize}
\end{itemize} 
\end{conj}

L'implication $(i) \Rightarrow (ii)$ de la Conjecture \ref{congetturaintro} découle de \cite[Lemma 7.9]{pas}. Donc il reste à montrer l'implication $(ii) \Rightarrow (i)$.
\medbreak
Le cas  $\alpha \in \Oe^{\times}$  (resp. $\beta \in \Oe^{\times}$) est facile (Proposition \ref{proposizionefacile}).  On peut donc supposer  $\alpha, \beta \notin \Oe^{\times}$. En reformulant $\rho_{\vec{d}}\otimes \pi$ comme induite compacte (Proposition \ref{indotte}), on en définit un sous-$\Oe[\GL(\OF)]$-module de type fini naturel $\Theta$ qui l'engendre sur $E$. Prouver la Conjecture \ref{congetturaintro} équivaut alors à prouver que le réseau $\Theta$ est séparé, i.e. il ne contient pas de $E$-droite. Dans cet article on montre que c'est bien le cas, une fois que certaines hypothèses supplémentaires sur le vecteur $\vec{d}$ sont satisfaites. 
\medbreak
Plus précisément, d'après la définition de $\Theta$ (voir la Section \ref{capo}) on dispose d'une application surjective de $\Oe[G]$-modules:
\[
\theta\colon \frac{\ind \rho_{\vec{d}}^0}{(T-a_p)(\ind \rho_{\vec{d}}^0)} \onto \Theta_{\vec{d},a_p}
\] 
où, $\rho_{\vec{d}}^0$ est un sous-$\Oe$-module de $\rho_{\vec{d}}$ qui l'engendre sur $E$, $a_p$ désigne un élément de $\mathfrak{p}_E$ qui dépend de $\Theta$ et $T$ est l'opérateur de Hecke usuel. Le résultat principal de l'article donne une condition nécessaire et suffisante sur le vecteur $\vec{d}$ pour que l'application $\theta$ soit injective. 

\begin{theo}\label{principe2}
Avec les notations précédentes l'application $\theta$ est injective (et donc un isomorphisme) si et seulement si les deux conditions suivantes sont satisfaites: 
\begin{itemize}
\item[(i)] Pour tout $l\in \{0,\ldots, f-1\}$ on a $|J_l|\leq 1$;
\item[(ii)] Si $\sigma \in J_l$ on a  
\[
d_{\sigma}+1 \leq p^{v_{\sigma}}.
\]
\end{itemize}
\end{theo}

En utilisant le Théorème \ref{principe2} on déduit le corollaire suivant.

\begin{cor}\label{libero2}
Supposons (1) et (2) ci-dessus, et supposons que $\vec{d} = (d_{\sigma})_{\sigma \in S}$ satisfait les conditions (i) et (ii) du Théorème \ref{principe2}. Alors, le $\Oe$-réseau $\Theta_{\vec{d},a_p}$ est séparé.
\end{cor}

\begin{rem}{\rm
Dans \cite{marco7} nous avons donné une description explicite du complété unitaire universel de la  représentation $\rho_{\vec{d}}\otimes \pi$ en utilisant certains espaces de Banach de fonctions de classe $C^r$ sur $\OF$ (\cite{marco}), où $r$ est un convenable nombre réel positif. Une conséquence immédiate du Corollaire \ref{libero2} est que ce complété unitaire universel est non nul. }   
\end{rem}


\addtocontents{toc}{\protect\setcounter{tocdepth}{2}}

  

\section{Préliminaires}
\subsection{Rappels sur l'arbre de Bruhat-Tits}
Nous renvoyons à \cite{serre} et \cite{breuilaa} pour plus de détails concernant la construction et les propriétés de l'arbre de Bruhat-Tits $\mathcal{T}$ de $\mathrm{SL}_2(F)$.

On note $V$ un $F$-espace vectoriel de dimension $2$ et l'on en fixe une base $(e_1,e_2)$.  L'espace $V$ est muni d'une action de $G$ définie pour tout $g = [\begin{smallmatrix} {a} & {b} \cr {c} & {d} \end{smallmatrix}] \in G$ par
\[
g  e_1 = a e_1 + c e_2, \quad g  e_2 = b e_1 + d e_2. 
\]
Si $f_1, f_2$ forment une base de $V$, on note $[f_1,f_2]$  la classe d'homothétie du $\OF$-réseau $\OF f_1 \oplus \OF f_2$. Le groupe $G$ agit transitivement sur l'ensemble $\mathcal{S}$ des classes d'homothéties de réseaux de $V$: si $[f_1,f_2]$ est une telle classe et $g \in G$, alors $g [f_1,f_2] = [g f_1,g f_2]$. Comme le stabilisateur de $[e_1,e_2]$ est $KZ$, l'ensemble des classes d'homothéties est isomorphe, en tant que $G$-ensemble, à $G/KZ$.

Soient $s,s' \in \mathcal{S}$. Si on choisit un réseau $\Lambda$ dans la classe d'homothétie $s$, il existe un représentant $\Lambda'$ de $s'$ et un seul tel que $\Lambda' \subset \Lambda$ et $\Lambda/\Lambda'$ est monogène. Pour un tel $\Lambda'$ on a
\[
\Lambda/\Lambda' \simeq \OF/\uni^n \OF.
\]  
L'entier $n$ ne dépend que des classes $s$ et $s'$ de $\Lambda$ et $\Lambda'$; on le note $d(s,s')$; c'est \textit{la distance de} $s$ à $s'$.

Soit $\mathcal{T}$ l'arbre de Bruhat-Tits de $\mathrm{SL}_2(F)$. Les sommets de $\mathcal{T}$ sont les classes d'homothétie de réseaux de $V$ et les arêtes orientées sont les paires $(s,s')$ avec $d(s,s') = 1$.

Rappelons que l'on a la décomposition de Cartan:
\begin{align*}
G = \coprod_{n \in \mathbb{Z}_{\geq 0}} KZ \alpha^{-n} KZ = \Big(\coprod_{n \in \mathbb{Z}_{\geq 0}} IZ \alpha^{-n}KZ \Big) \coprod \Big(\coprod_{n \in \mathbb{Z}_{\geq 0}} IZ \beta \alpha^{-n}KZ \Big).
\end{align*}
Alors, pour tout $n \in \mathbb{Z}_{\geq 0}$, les classes de $KZ \alpha^{-n} KZ/KZ$ correspondent aux sommets $s_i$ de $\mathcal{T}$ tels que $d(s_i, s_0) = n$ où $s_0 \ugu [e_1,e_2]$. Posons $I_0 = \{0 \}$ et si $n \in \mathbb{Z}_{>0}$
\[
I_n = \big\{[\mu_0]+ \uni [\mu_1]+\ldots +\uni^{n-1} [\mu_{n-1}], \ (\mu_0,\ldots,\mu_{n-1}) \in (k_F)^n   \big\} \subseteq \OF.
\] 
Pour $n \in \mathbb{Z}_{\geq 0}$ et $\mu \in I_n$ posons:
\begin{align*}
g_{n,\mu}^0 = \begin{bmatrix} {\uni^n} & {\mu} \cr {0} & {1} \end{bmatrix}, \quad 
g_{n,\mu}^1 = \begin{bmatrix} {1} & {0} \cr {\uni \mu} & {\uni^{n+1}} \end{bmatrix}.
\end{align*}
En particulier, $g_{0,0}^0$ est la matrice identité, $g_{0,0}^1 = \alpha$ et, pour tout $n \in \mathbb{Z}_{\geq 0}$, $\mu \in I_n$ on a $g_{n,\mu}^1 = \beta g_{n,\mu}^0 w$. Par ailleurs, les $g_{n,\mu}^0$ et $g_{n,\mu}^1$ forment un système de représentants de $G/KZ$:
\begin{align}\label{cartan}
G = \Big(\coprod_{n \in \mathbb{Z}_{\geq 0}, \mu \in I_n} g_{n,\mu}^0 KZ \Big) \coprod \Big( \coprod_{n \in \mathbb{Z}_{\geq 0}, \mu \in I_n} g_{n,\mu}^1 KZ \Big); 
\end{align}
Pour $n\in \mathbb{Z}_{\geq 0}$ posons:
\begin{align*}
S_n^0 = IZ \alpha^{-n}KZ = \coprod_{\mu \in I_n} g_{n,\mu}^0 KZ, \quad
S_n^1 = IZ \beta\alpha^{-n}KZ = \coprod_{\mu \in I_n} g_{n,\mu}^1 KZ 
\end{align*}
comme dans \cite{breuilaa}. Notons $S_n = S_n^0 \coprod S_n^1$ et $B_n = B_n^0 \coprod B_n^1$, où $B_n^0 = \coprod_{m \leq n} S_m^0$ et  $B_n^1 = \coprod_{m \leq n} S_m^1$. En particulier $S_0 = KZ \coprod \alpha KZ$.

\begin{rem}[\cite{breuilaa}]
{\rm On peut voir que $S_n^0 \coprod S_{n-1}^1$ (resp. $B_n^0 \coprod B_{n-1}^1$) est l'ensemble des sommets de $\mathcal{T}$ de distance $n$ (resp. inférieure ou égale à $n$) de $KZ$. De même $S_n^1 \coprod S_{n-1}^0$ (resp. $B_n^1 \coprod B_{n-1}^0$) est l'ensemble des sommets de $\mathcal{T}$ de distance $n$ (resp. inférieure ou égale à $n$) de $\alpha KZ$.} 
\end{rem}

On désigne par $R$ le corps $E$ ou son anneau d'entiers $\Oe$. Soit $\sigma$ une représentation $R$-linéaire de $KZ$ sur un $R$-module libre $V_{\sigma}$ de rang fini. On note:
\[
\ind \sigma
\]
la représentation de $G$ sur $R$ dont le $R$-module sous-jacent est l'ensemble des fonctions $f\colon G \to V_{\sigma}$ à support compact modulo $Z$ telles que 
\[
\forall \kappa \in KZ, \forall g \in G, \quad f(\kappa g) = \sigma(\kappa)f(g)
\] 
et sur lesquelles $G$ agit par translation à droite (i.e. $(g\cdot f)(g') \ugu f(g'g)$). Comme dans \cite{breuilab}, pour $g \in G$ et $v \in V_{\sigma}$ on note $[g,v]$ l'élément de $\ind\sigma$ défini comme suit:
\begin{align*}
[g,v](g') = \left\{ \begin{array}{ll}
 \sigma(g'g)(v)   & \mbox{si} \ g' \in KZg^{-1} \\
0 & \mbox{si}\ g' \notin KZg^{-1} 
\end{array}\right.
\end{align*}
On a les égalités suivantes $g_1  [g_2, v] = [g_1 g_2, v]$ et $[g\kappa, v] = [g,\sigma(\kappa)(v)]$ si $g_1,g_2,g \in G$ et $\kappa \in KZ$.

Rappelons le résultat suivant qui décrive une base de la représentation $\ind \sigma$.
\begin{prop}\label{baseind}
Soit $\mathcal{B}$ une $R$-base de $V_{\sigma}$ et $\mathcal{G}$ un système de représentants des classe à gauche de $G/KZ$. Alors la famille de fonctions:
\[
\mathcal{I} \ugu \{[g,v], \ g\in \mathcal{G}, v\in \mathcal{B} \}
\]
forme une base de $\ind\sigma$.
\begin{proof}
Voir \cite[§2]{bl}.
\end{proof}
\end{prop}

\begin{rem}\label{prodottotens}
{\rm La représentation $\ind\sigma$ est isomorphe à la représentation de $G$ portée par le $R[G]$-module $R[G]\otimes_{R[KZ]}V_{\sigma}$. Plus précisément, si $g \in G$ et $v \in V_{\sigma}$, alors l'élément $g \otimes v$ correspond à la fonction $[g,v]$.} 
\end{rem}

D'après la Proposition \ref{baseind} et d'après \eqref{cartan} toute fonction $f\in \ind\sigma$ s'écrit de façon unique comme somme finie:
\[
f = \sum_{n = 0}^{\infty}\sum_{\mu \in I_n}  \big([g_{n,\mu}^0,v_{n,\mu}^0]+ [g_{n,\mu}^1,v_{n,\mu}^1]\big), 
\] 
où $v_{n,\mu}^0$, $v_{n,\mu}^1 \in V_{\sigma}$. On appelle support de $f$ l'ensemble des  $g_{n,\mu}^i$ tels que $v_{n,\mu}^i \neq 0$. On écrit $f \in S_n$ si le support de $f$ est contenu dans $S_n$ et de manière analogue pour $B_n, S_n^0$, etc.

Soit $\pi$ une représentation $R$-linéaire de $G$ sur un $R$-module. Par la Remarque \ref{prodottotens} et par \cite[Théorème 2.19]{curtis} on déduit un isomorhisme canonique de $R$-modules:
\[
\mathrm{Hom}_{R[G]}(\ind\sigma, \pi) \simeq \mathrm{Hom}_{R[H]}(\sigma,\pi|_{KZ}).
\]
Autrement dit, le foncteur d'induction $\ind$ est un adjoint à gauche du foncteur de restriction. Cet énoncé est connu sous le nom de \textit{réciprocité de Frobenius}. 

\subsection{Rappels sur les algèbres de Hecke}

Soit $\sigma$ une représentation $R$-linéaire de $KZ$ sur un $R$-module libre $V_{\sigma}$ de rang fini. L'algèbre de Hecke $\mathcal{H}(KZ,\sigma)$ associée à $KZ$ et $\sigma$ est définie par:
\[
\mathcal{H}(KZ,\sigma) = \mathrm{End}_{R[G]}(\ind\sigma).
\] 
C'est une $R$-algèbre;

On peut réinterpréter $\mathcal{H}(KZ,\sigma)$ comme une algèbre de convolution: notons $\mathcal{H}_{KZ}(\sigma)$ le $R$-module des fonctions $\varphi\colon G \longrightarrow \mathrm{End}_R(V_{\sigma})$ à support compact modulo $Z$ et telles que
\begin{align*}
\forall \kappa_1,\kappa_2 \in KZ, \forall g \in G, \quad    \varphi(\kappa_1 g \kappa_2) = \sigma(\kappa_1) \circ \varphi(g) \circ \sigma(\kappa_2).
\end{align*}
C'est une $R$-algèbre unifère pour le produit de convolution défini pour tout $\varphi_1,\varphi_2 \in \mathcal{H}_{KZ}(\sigma)$ et pour tout $g \in G$ par la formule
\[
\varphi_1 * \varphi_2(g) = \sum_{xKZ \in G/KZ} \varphi_1(x) \varphi_2(x^{-1}g),
\]
et d'élément unité la fonction
\begin{align*}
\varphi_e(g) = \left\{ \begin{array}{ll}
\sigma(g)   & \mbox{si} \ g \in KZ, \\
0 & \mbox{sinon}.
\end{array} \right. 
\end{align*}
On vérifie facilement que l'accouplement
\begin{align*}
\mathcal{H}_{KZ}(\sigma) \times \ind\sigma &\longrightarrow \ind\sigma \\
(\varphi,f) &\longmapsto \left\langle \varphi,f \right\rangle(g) \ugu \sum_{xKZ \in G/KZ} \varphi(x)(f(x^{-1}g)),
\end{align*}
munit $\ind\sigma$ d'une structure de $\mathcal{H}_{KZ}(\sigma)$-module à gauche qui commute à l'action de $G$.

\begin{lemma}\label{isoalg}
L'application:
\begin{align*}
\mathcal{H}_{KZ}(\sigma) &\longrightarrow \mathcal{H}(KZ,\sigma) \\
\varphi &\longmapsto T_{\varphi}(f) \ugu \left\langle \varphi,f \right\rangle
\end{align*}
est un isomorphisme de $R$-algèbres. En particulier, si $g \in G$ et $v \in V_{\sigma}$  l'action de $T_{\varphi}$ sur $[g,v]$ est donnée par
\begin{align}\label{hecke}
T_{\varphi}([g,v]) = \sum_{xKZ \in G/KZ} [gx,\varphi(x^{-1})(v)].
\end{align}
\begin{proof}
Ce résultat découle de \cite[Proposition 5]{bl} dans le cas où $\sigma$ est une représentation lisse. Le cas général résulte de \cite[Lemme 1.2]{st}. 
\end{proof}
\end{lemma}

Supposons $R=E$. On désigne par $\mathbf{1}$ la représentation triviale de $KZ$ et on suppose que $\sigma$ est la restriction à $KZ$ d'une représentation localement $\Q$-analytique de $G$ (au sens de \cite{sch1}, \cite{st2}) sur $V_{\sigma}$. L'application:
\begin{align*}
\iota_{\sigma}\colon \he_{KZ}(\mathbf{1}) &\longrightarrow \he_{KZ}(\sigma) \\
\varphi &\longmapsto (\varphi\cdot \sigma)(g) \ugu \varphi(g)\sigma(g)
\end{align*}
est un homomorphisme injectif de $E$-algèbres (\cite[§1]{st}). Le résultat suivant donne une condition suffisante pour que l'application $\iota_{\sigma}$ soit bijective. 

Rappelons d'abord que l'on a une action $\Q$-linéaire de l'algèbre de Lie $\mathfrak{g}$ de $G$ sur l'espace $V_{\sigma}$ définie par:
\[
\mathfrak{x} v = \frac{d}{dt}\exp (t\mathfrak{x}) v|_{t=0}
\]
où $\exp\colon \mathfrak{g} \dashrightarrow G$ désigne l'application exponentielle définie localement autour de $0$ (\cite[§2]{sch1}). Cette action se prolonge en une action de l'algèbre de Lie $\mathfrak{g}\otimes_{\Q}E$. 
\begin{lemma}\label{liealg} Supposons que le $\mathfrak{g}\otimes_{\Q}E$-module $V_{\sigma}$ soit absolument irréductible; alors l'application $\iota_{\sigma}$ est bijective.
\begin{proof}
Ce résultat est démontré dans \cite[Lemme 1.4]{st} dans le cas particulier où $V_{\sigma}$ est une représentation $F$-analytique et $F$ est un sous-corps de $E$. Les mêmes arguments s'appliquent \textit{mutatis mutandis}.
\end{proof}
\end{lemma}

\section{Représentations de $\GL(F)$}

\subsection{Représentations $\Q$-algébriques de $\GL(F)$}\label{repr}
On note  $(\mathrm{Sym}^{d_{\sigma}}E^2)^{\sigma}$, pour $\sigma \in S$ et $d_{\sigma} \in \mathbb{Z}_{\geq 0}$, la représentation algébrique irréductible de $\mathrm{GL}_2 \otimes_{F,\sigma} E$ dont le plus haut poids est $\chi_{\sigma}\colon \mathrm{diag}(x_1,x_2) \mapsto \sigma(x_2)^{d_{\sigma}}$ vis-à-vis du sous-groupe des matrices triangulaires supérieures.  Si $d_{\sigma}$ est impair, on choisit une racine carrée de $\sigma(\uni)$. Notons $\chi\colon \GL(F) \to F^{\times}$ le caractère défini par:
\[
\begin{bmatrix} {a} & {b} \cr {c} & {d} \end{bmatrix} \mapsto \uni^{- val_F\big(\mathrm{d\acute{e}t}\big([\begin{smallmatrix} {a} & {b} \cr {c} & {d} \end{smallmatrix}]\big)\big)/f}
\]
et posons:
\[
\big(\Sym^{d_{\sigma}} E^2\big)^{\sigma} = (\Sym^{d_{\sigma}} E^2)^{\sigma} \otimes_E \big(\sigma \circ \chi\big)^{\frac{d_{\sigma}}{2}}.
\]  
On identifie $\big(\Sym^{d_{\sigma}} E^2\big)^{\sigma}$ avec le $E$-espace vectoriel $\bigoplus_{i_{\sigma} = 0}^{d_{\sigma}} E x_{\sigma}^{d_{\sigma}- i_{\sigma}} y_{\sigma}^{i_{\sigma}}$ des polynômes homogènes de degré $d_{\sigma}$ en $x_{\sigma}$ et $y_{\sigma}$ à coefficients dans $E$ avec l'action à gauche de $G$ donnée explicitement pour tout $0\leq i_{\sigma} \leq d_{\sigma}$ par
\begin{align}\label{azione}
\begin{bmatrix} {a} & {b} \cr {c} & {d} \end{bmatrix}  (x_{\sigma}^{d_{\sigma}-i_{\sigma}} y_{\sigma}^{i_{\sigma}}) = \big(  \sigma\circ \chi([\begin{smallmatrix} {a} & {b} \cr {c} & {d} \end{smallmatrix}])\big)^{\frac{d_{\sigma}}{2}} (\sigma(a)x_{\sigma}+ \sigma(c)y_{\sigma})^{d_{\sigma}-i_{\sigma}} (\sigma(b)x_{\sigma}+ \sigma(d)y_{\sigma})^{i_{\sigma}}.
\end{align}
Si $v_{\sigma} \in \big(\Sym^{d_{\sigma}} E^2\big)^{\sigma}$ et $g\in G$, on note simplement $g v_{\sigma}$ l'action de $g$ sur $v_{\sigma}$. 

\begin{rem}
{\rm Par la formule \eqref{azione} on déduit:
\[
\begin{bmatrix} {\uni} & {0} \cr {0} & {\uni} \end{bmatrix}  v_{\sigma} = v_{\sigma}  
\]
pour tout $v_{\sigma} \in \big(\Sym^{d_{\sigma}} E^2\big)^{\sigma}$.
} 
\end{rem}

Choisissons une numérotation $\sigma_1,\ldots,\sigma_{|S|}$ des plongements de $F$ dans $E$. On fixe $\vec{d} = (d_{\sigma_1},\ldots, d_{\sigma_{|S|}} )$ un $|S|$-uplet d'entiers positifs ou nuls. Posons: 
\[
I_{\vec{d}} = \big\{(i_{\sigma_1},\ldots, i_{\sigma_{|S|}}) \in \mathbb{Z}_{\geq 0}^{|S|},\,  0\leq i_{\sigma_j} \leq d_{\sigma_{j}} \ \mbox{pour tout} \ 1 \leq j \leq |S| \big\}
\]
et munissons-le de l'ordre lexicographique que l'on notera $\prec$.

Notons $\rho_{\vec{d}}$ l'unique  représentation de $G$ dont l'espace sous-jacent est 
\[
V_{\rho_{\vec{d}}} = \big(\Sym^{d_{\sigma_1}} E^2\big)^{\sigma_1}  \otimes_E \big(\Sym^{d_{\sigma_2}}E^2\big)^{\sigma_2} \otimes_E \ldots \otimes_E \big(\Sym^{d_{\sigma_{|S|}}}E^2\big)^{\sigma_{|S|}}
\]
et sur lequel un élément $[\begin{smallmatrix} {a} & {b} \cr {c} & {d} \end{smallmatrix}]$ de $G$ agit par
\begin{align}\label{azione2}
\rho_{\vec{d}} \Big(\begin{bmatrix} {a} & {b} \cr {c} & {d} \end{bmatrix}\Big) (v_{\sigma_1}\otimes \ldots \otimes v_{\sigma_{|S|}}) = \begin{bmatrix} {a} & {b} \cr {c} & {d} \end{bmatrix}   v_{\sigma_1}\otimes \begin{bmatrix} {a} & {b} \cr {c} & {d} \end{bmatrix}   v_{\sigma_2}\otimes  \ldots \otimes \begin{bmatrix} {a} & {b} \cr {c} & {d} \end{bmatrix}   v_{\sigma_{|S|}}.
\end{align}

\begin{rem}\label{absirr}
{\rm La représentation $\rho_{\vec{d}}$ est une représentation absolument irréductible de $G$ qui reste absolument irréductible sous l'action restreinte d'un sous-groupe ouvert de $G$ (voir \cite[§2]{bs}).}
\end{rem}

Posons pour tout $\vec{i}\in I_{\vec{d}}$:
\[
e_{\vec{d},{\vec{i}}} = e_{d_{\sigma_1},i_{\sigma_1}} \otimes \ldots \otimes e_{d_{\sigma_{|S|}},i_{\sigma_{|S|}}},
\]
où, pour tout $1\leq j\leq |S|$, $e_{d_{\sigma_j},i_{\sigma_j}}$ désigne le monôme  $x_{\sigma_j}^{d_{\sigma_j}-i_{\sigma_j}} y_{\sigma_j}^{i_{\sigma_j}}$, et notons  $U_{\vec{d}}$ l'endomorphisme  de $V_{\rho_{\vec{d}}}$ défini par
\[
U_{\vec{d}} = U_{d_{\sigma_1}} \otimes \ldots \otimes  U_{d_{\sigma_{|S|}}}, 
\]
où, pour tout $\sigma \in S$, $U_{d_{\sigma}}$ désigne l'endomorphisme de  $(\Sym^{d_{\sigma}} E^2\big)^{\sigma}$ explicitement donné dans la base $(e_{d_{\sigma},i_{\sigma}})_{0\leq i_{\sigma} \leq d_{\sigma}}$ par la matrice diagonale
\begin{align}\label{matricediag}
U_{d_{\sigma}} =
\begin{bmatrix}
\sigma(\uni)^{d_{\sigma}} & 0 & \ldots & 0 \\
0 & \sigma(\uni)^{d_{\sigma}-1} & \ddots & \vdots \\
\vdots &  \ddots &  \ddots & 0 \\
0 &  \ldots &  0 & 1
\end{bmatrix}
\end{align}

Choisissons pour tout $\sigma \in S$ tel que $d_{\sigma}$ est impair une racine carrée de $\sigma(\uni)$ dans $E$. La notation $\sigma(\uni)^{\frac{a_{\sigma}}{2}}$ a donc un sens pour $a_{\sigma} \in \mathbb{Z}$.

\begin{lemma}\label{funzione}
Il existe une unique fonction $\psi \colon G \longrightarrow \mathrm{End}_{E}(V_{\rho_{\vec{d}}})$ à support dans $KZ \alpha^{-1} KZ$ telle que
\begin{itemize}
\item[(i)] pour tout $\kappa_1,\kappa_2 \in KZ$ on a $\psi(\kappa_1 \alpha^{-1}\kappa_2) = {\rho_{\vec{d}}}(\kappa_1) \circ \psi(\alpha^{-1}) \circ \rho_{\vec{d}}(\kappa_2)$;
\item[(ii)] dans la base $\{e_{\vec{d},\vec{i}}, \vec{i} \in I_{\vec{d}} \}$, $\psi(\alpha^{-1}) = U_{\vec{d}}$.    
\end{itemize}
\begin{proof}
\underline{\textit{Existence}}. Il faut vérifier que $\psi$ est bien définie, c'est-à-dire si $\kappa_1\alpha^{-1}\kappa_2 = \alpha^{-1}$ dans $G$ pour $\kappa_1,\kappa_2 \in KZ$, alors $\psi(\kappa_1\alpha^{-1}\kappa_2) = \psi(\alpha^{-1})$. Notons $u_{\sigma} = \frac{d_{\sigma}}{2}$ pour tout $\sigma \in S$. D'après les formules \eqref{azione2} et \eqref{azione} on a pour tout $\vec{i} \in I_{\vec{d}}$
\[
\rho_{\vec{d}}(\alpha^{-1})(e_{\vec{d},\vec{i}}) = \uni^{\vec{u} - \vec{i}} e_{\vec{d},\vec{i}}
\]
et donc on déduit que l'automorphisme $\rho_{\vec{d}}(\alpha^{-1})$ de l'espace $V_{\rho_{\vec{d}}}$ est donné, dans la base $\{e_{\vec{d},\vec{i}}, \vec{i} \in I_{\vec{d}} \}$,  par la matrice $\uni^{-\vec{u}}\psi(\alpha^{-1})$. Donc, si  $\kappa_1\alpha^{-1}\kappa_2 = \alpha^{-1}$ dans $G$ pour $\kappa_1,\kappa_2 \in KZ$, on  déduit:
\[
\rho_{\vec{d}}(\kappa_1)\circ \uni^{-\vec{u}} \psi(\alpha^{-1})  \circ \rho_{\vec{d}}(\kappa_2) = \uni^{-\vec{u}}\psi(\alpha^{-1}),
\]
d'où 
\[
\rho_{\vec{d}}(\kappa_1)\circ \psi(\alpha^{-1})  \circ \rho_{\vec{d}}(\kappa_2) = \psi(\kappa_1 \alpha^{-1} \kappa_2) = \psi(\alpha^{-1}).
\]

\underline{\textit{Unicité}}. Claire.
\end{proof}
\end{lemma}

D'après le Lemme \ref{isoalg} on sait que l'algèbre de Hecke $\mathcal{H}(KZ,\rho_{\vec{d}})$ est naturellement isomorphe à l'algèbre de convolution $\mathcal{H}_{KZ}(\rho_{\vec{d}})$ des fonctions $\varphi\colon G \longrightarrow \mathrm{End}_E(V_{\rho_{\vec{d}}})$ à support compact modulo $Z$ et telles que
\begin{align*}
\forall \kappa_1,\kappa_2 \in KZ, \forall g \in G, \quad \varphi(\kappa_1 g \kappa_2) = \rho_{\vec{d}}(\kappa_1) \circ \varphi(g) \circ \rho_{\vec{d}}(\kappa_2).
\end{align*}
À l'application $\psi$ du Lemme \ref{funzione} est donc en particulier associé un opérateur $T \in \mathcal{H}(KZ,\rho_{\vec{d}})$.  L'action de $T$ sur les éléments $[g,v]$, pour $g \in G$ et $v \in V_{\rho_{\vec{d}}}$, est donnée par la formule \eqref{hecke}.

On peut énoncer le résultat bien connu suivant. 

\begin{lemma}\label{ophec}
Il existe un isomorphisme de $E$-algèbres:
\[
\mathcal{H}(KZ,\rho_{\vec{d}}) \simeq E[T].
\]
\begin{proof}
Comme l'espace $V_{\rho_{\vec{d}}}$ est un $\mathfrak{gl}_2(F) \otimes_{\Q}E$-module absolument irréductible, où $\mathfrak{gl}_2(F)$ désigne l'algèbre de Lie de $\GL(F)$, le Lemme \ref{liealg} implique que l'application
\begin{align*}
\iota_{\rho_{\vec{d}}}\colon \he_{KZ}(\mathbf{1}) &\longrightarrow \he_{KZ}(\rho_{\vec{d}}) \\
f &\longmapsto (f\cdot \rho_{\vec{d}})(g) \ugu f(g)\rho_{\vec{d}}(g)
\end{align*} 
est un isomorphisme de $E$-algèbres. En utilisant le Lemme \ref{isoalg} on déduit qu'il existe une unique application de $E$-algèbres $u_{\rho_{\vec{d}}} \colon \he(KZ,\mathbf{1}) \to \he(KZ,\rho_{\vec{d}})$  de sorte que le diagramme suivant
\[
\xymatrix{
\he_{KZ}(\mathbf{1}) \ar[r]^{\sim} \ar[d]_{\iota_{\rho_{\vec{d}}}}
& \he(KZ,\mathbf{1})  \ar[d]_{u_{\rho_{\vec{d}}}} \\
 \he_{KZ}(\rho_{\vec{d}}) \ar[r]^{\sim} & \he(KZ,\rho_{\vec{d}})  }
\]
soit commutatif. Par construction,  l'application $u_{\rho_{\vec{d}}}$ est un isomorphisme de $E$-algèbres. Notons $T_1$ l'élément de  $\he(KZ,\mathbf{1})$ qui correspond à $\mathbf{1}_{KZ\alpha^{-1}KZ}$ via la réciprocité de Frobenius. Un raisonnement analogue à celui de \cite[Proposition 4]{bl2} montre que  $\he(KZ,\mathbf{1})$ est isomorphe à l'algèbre de polynômes $E[T_1]$. Or $u_{\rho_{\vec{d}}}(T_1) = (\prod_{\sigma \in S}\sigma(\uni)^{-\frac{d_{\sigma}}{2}}) T$, d'où le résultat. 
\end{proof}
\end{lemma}

\begin{rem}\label{iniet}
{\rm Un raisonnement sans difficultés sur l'arbre de Bruhat-Tits de $\mathrm{SL}_2(F)$ montre que l'opérateur $T$ est injectif sur l'espace $\ind \rho_{\vec{d}}$.  }
\end{rem}

\subsection{Réseaux}\label{lattici}
Conservons les notations de la Section \ref{repr} et notons  $(\underline{\Sym}^{d_{\sigma}} \Oe^2)^{\sigma}$, pour $\sigma \in S$ et $d_{\sigma} \in \mathbb{Z}_{\geq 0}$, la représentation du groupe $KZ$ ayant pour espace sous-jacent le $\Oe$-module $\bigoplus_{i_{\sigma} = 0}^{d_{\sigma}} \Oe x_{\sigma}^{d_{\sigma}- i_{\sigma}} y_{\sigma}^{i_{\sigma}}$ des polynômes homogènes de degré $d_{\sigma}$, sur lequel un élément $[\begin{smallmatrix} {a} & {b} \cr {c} & {d} \end{smallmatrix}] \in K$ agit par
\begin{align}\label{azione3}
\begin{bmatrix} {a} & {b} \cr {c} & {d} \end{bmatrix} (x_{\sigma}^{d_{\sigma}-i_{\sigma}} y_{\sigma}^{i_{\sigma}}) =  (\sigma(a)x_{\sigma}+ \sigma(c)y_{\sigma})^{d_{\sigma}-i_{\sigma}} (\sigma(b)x_{\sigma}+ \sigma(d)y_{\sigma})^{i_{\sigma}}
\end{align}
et la matrice $[\begin{smallmatrix} {\uni} & {0} \cr {0} & {\uni} \end{smallmatrix}] \in Z$ agit comme l'identité. Si $v_{\sigma} \in \big(\Sym^{d_{\sigma}} \Oe^2\big)^{\sigma}$ et $g\in G$, on note simplement $g v_{\sigma}$ l'action de $g$ sur $v_{\sigma}$.

\begin{defin}\label{lattice}
Soit $V$ un $E$-espace vectoriel. Un \textit{réseau} $\mathcal{L}$ de $V$ est un sous-$\Oe$-module de $V$ tel que pour tout $v \in V$ il existe un élément non nul $a \in E^{\times}$ tel que $av \in \mathcal{L}$. On dit qu'un réseau $\mathcal{L}$ est \textit{séparé} si $\bigcap_{n \in \mathbb{N}} \varpi_{E}^n \mathcal{L} = 0$ ou, de manière équivalente, s'il ne contient pas de $E$-droite.
\end{defin}

\begin{rem}\label{sottorep}
{\rm Le $\Oe$-module $(\Sym^{d_{\sigma}} \Oe^2)^{\sigma}$ est un réseau séparé de  $(\Sym^{d_{\sigma}} E^2)^{\sigma}$ qui est  stable sous l'action de $KZ$. }
\end{rem}

Notons $\rho_{\vec{d}}^0$ l'unique représentation de $KZ$ dont l'espace sous-jacent est défini par 
\[
V_{\rho_{\vec{d}}^0} = \big(\Sym^{d_{\sigma_1}} \Oe^2\big)^{\sigma_1}  \otimes_{\Oe} \big(\Sym^{d_{\sigma_2}}\Oe^2\big)^{\sigma_2} \otimes_{\Oe} \ldots \otimes_{\Oe} \big(\Sym^{d_{\sigma_{|S|}}}\Oe^2\big)^{\sigma_{|S|}},
\]
et sur lequel un élément $[\begin{smallmatrix} {a} & {b} \cr {c} & {d} \end{smallmatrix}]$ de $KZ$ agit par
\begin{align}\label{azione4}
\rho_{\vec{d}}^0 \Big(\begin{bmatrix} {a} & {b} \cr {c} & {d} \end{bmatrix}\Big) (v_{\sigma_1}\otimes \ldots \otimes v_{\sigma_{|S|}}) = \begin{bmatrix} {a} & {b} \cr {c} & {d} \end{bmatrix}  v_{\sigma_1}\otimes \begin{bmatrix} {a} & {b} \cr {c} & {d} \end{bmatrix}   v_{\sigma_2}\otimes  \ldots \otimes \begin{bmatrix} {a} & {b} \cr {c} & {d} \end{bmatrix}   v_{\sigma_{|S|}}.
\end{align}
La Remarque \ref{sottorep} implique que le $\Oe$-module $V_{\rho_{\vec{d}}^0}$ est un réseau séparé de l'espace $V_{\rho_{\vec{d}}}$ construit dans le paragraphe \ref{repr} qui est  stable sous l'action de $KZ$. Donc le $\Oe$-module $\ind \rho_{\vec{d}}^0$ est encore un réseau séparé de l'espace $\ind \rho_{\vec{d}}$, stable sous l'action de $G$.

Par extension des scalaires on déduit alors une application injective de $\he(KZ,\rho_{\vec{d}}^0)$ dans  $\he(KZ,\rho_{\vec{d}})$. De plus, l'opérateur $T \in \he(KZ,\rho_{\vec{d}})$ décrit dans le paragraphe \ref{repr} induit par restriction un endomorphisme $G$-équivariant de $\ind{\rho_{\vec{d}}^0}$ que l'on note encore $T$. D'après le Lemme \ref{ophec} on déduit alors un isomorphisme de $\Oe$-algèbres entre  $\he(KZ,\rho_{\vec{d}}^0)$ et l'algèbre de polynômes $\Oe[T]$. En utilisant le Lemme \ref{isoalg}, ce qui précède peut se résumer dans le diagramme commutatif suivant:
\[
\xymatrix{
\he_{KZ}(\rho_{\vec{d}}^0) \ar[r]^{\sim} \ar@{^{(}->}[d] & \he(KZ,\rho_{\vec{d}}^0) \ar@{^{(}->}[d] & \hspace{-27pt}\simeq \Oe[T]  \\
 \he_{KZ}(\rho_{\vec{d}}) \ar[r]^{\sim} & \he(KZ,\rho_{\vec{d}}) & \hspace{-33pt} \simeq E[T] }
\]


\subsection{Formulaire}\label{formulario}
Fixons $\vec{d}$ un $|S|$-uplet d'entiers positifs ou nuls et reprenons les notations des Sections \ref{repr} et \ref{lattici}. Pour $0\leq m \leq n$, soit $[\ ]_m \colon I_n \to I_m$ les applications ``troncatures" définies par:
\begin{align*}
\Big[\sum_{i=0}^{n-1} \uni^i [\mu_i] \Big]_m &= \left\{ \begin{array}{ll}
 \sum_{i=0}^{m-1} \uni^i [\mu_i]  & \mbox{si} \ m \geq 1, \\
0 & \mbox{si} \ m = 0.
\end{array} \right. 
\end{align*}
Rappelons que $\psi$ désigne la fonction définie dans le Lemme \ref{funzione}.

\begin{lemma}\label{tecnico1}
Soit $n \in \mathbb{Z}_{\geq 0}$, $\mu \in I_n$ et $v \in V_{\rho_{\vec{d}}^0}$. On a:
\[
T([g_{n,\mu}^0,v]) = T^+([g_{n,\mu}^0,v]) + T^-([g_{n,\mu}^0,v]), 
\]
où 
\begin{align*}
T^+([g_{n,\mu}^0,v]) &= \sum_{\lambda \in I_1} [g_{n+1,\mu + \uni^n \lambda }^0,(\rho_{\vec{d}}^0(w) \circ \psi(\alpha^{-1})\circ \rho_{\vec{d}}^0(w_{\lambda}))(v)], \\
T^-([g_{n,\mu}^0,v])                  &= \left\{ \begin{array}{ll}
 [g_{n-1,[\mu]_{n-1} }^0,(\rho_{\vec{d}}^0(w  w_{([\mu]_{n-1}-\mu)/\uni^{n-1}} ) \circ \psi(\alpha^{-1}))(v) ]    & \mbox{si} \ n \geq 1, \\
 \big[\alpha ,\psi(\alpha^{-1})(v)\big] & \mbox{si} \ n = 0.
\end{array} \right. 
\end{align*}
\begin{proof}
Voir \cite[§2.5]{breuilaa}.
\end{proof}
\end{lemma}

\begin{lemma}\label{tecnico2}
Soit $n \in \mathbb{Z}_{\geq 0}$, $\mu \in I_n$ et $v \in V_{\rho_{\vec{d}}^0}$. On a:
\[
T([g_{n,\mu}^1,v]) = T^+([g_{n,\mu}^1,v]) + T^-([g_{n,\mu}^1,v]), 
\]
où 
\begin{align*}
T^+([g_{n,\mu}^1,v]) &= \sum_{\lambda \in I_1} [g_{n+1,\mu + \uni^n \lambda }^1,(\psi(\alpha^{-1}) \circ \rho_{\vec{d}}^0(w_{\lambda}w))(v)], \\
T^-([g_{n,\mu}^1,v])                  &= \left\{ \begin{array}{ll}
 [g_{n-1,[\mu]_{n-1} }^0,(\rho_{\vec{d}}^0(  w_{([\mu]_{n-1}-\mu)/\uni^{n-1}}  ) \circ \psi(\alpha^{-1}) \circ \rho_{\vec{d}}^0(w))(v)]                   & \mbox{si} \ n \geq 1, \\
\phantom{} [ \mathrm{Id} ,(\rho_{\vec{d}}^0(w) \circ \psi(\alpha^{-1}) \circ \rho_{\vec{d}}^0(w))(v)] & \mbox{si} \ n = 0.
\end{array} \right. 
\end{align*}
\begin{proof}
Voir \cite[§2.5]{breuilaa}.
\end{proof}
\end{lemma}

D'après les Lemmes \ref{tecnico1} et \ref{tecnico2} on déduit facilement les deux égalités suivantes:
\begin{align*}
T^+([g_{n,\mu}^1,v]) &= \beta T^+([ g_{n,\mu}^0, \rho_{\vec{d}}^0(w)(v) ]), \\
T^-([g_{n,\mu}^1,v]) &= \beta  T^-([ g_{n,\mu}^0, \rho_{\vec{d}}^0(w)(v) ]).
\end{align*}

\begin{cor}\label{tecnico4}
Soit $n \in \mathbb{Z}_{\geq 0}$, $\mu,\lambda \in I_n$, $i,j \in \{0,1 \}$ et $v_1,v_2 \in V_{\rho_{\vec{d}}^0}$. Si $i \neq j$ ou $\mu \neq \lambda$ alors le support de $T^+([g_{n,\mu}^i,v_1])$ et le support de $T^+([g_{n,\lambda}^j,v_2])$ sont disjoints.
\begin{proof}
C'est une conséquence immédiate des lemmes \ref{tecnico1} et \ref{tecnico2}.
\end{proof}
\end{cor}

\begin{lemma}\label{tecnico5}
Soit $v = \sum_{\vec{0}\leq \vec{i}\leq \vec{d}} c_{\vec{i}} e_{\vec{d},\vec{i}} \in V_{\rho_{\vec{d}}^0}$ et $\lambda \in \OF$. On a:
\begin{align*}
&(\rho_{\vec{d}}(w) \circ \psi(\alpha^{-1}) \circ \rho_{\vec{d}}(w_{\lambda}))(v) = \sum_{\vec{0}\leq \vec{j}\leq \vec{d}} \Biggl(\uni^{\vec{j}} \sum_{\vec{j}\leq \vec{i}\leq \vec{d}} c_{\vec{i}} \binom{\vec{i}}{\vec{j}} (-\lambda)^{\vec{i}-\vec{j}} \Biggr) e_{\vec{d},\vec{j}}.
\end{align*}
\begin{proof}
D'après la formule \eqref{azione3} on a pour tout $\sigma \in S$ et tout $0\leq i_{\sigma} \leq d_{\sigma}$: 
\begin{align}\label{condizione171}
(w \circ U_{d_{\sigma}} \circ w_{\lambda})(e_{d_{\sigma},i_{\sigma}}) = \sum_{j_{\sigma} = 0}^{i_{\sigma}} \binom{i_{\sigma}}{j_{\sigma}} \sigma(\uni)^{j_{\sigma}}\sigma(-\lambda)^{i_{\sigma}-j_{\sigma}} e_{d_{\sigma},j_{\sigma}}.
\end{align}
Par \eqref{condizione171} et d'après la formule \eqref{azione4} on a:
\begin{align*}
(\rho_{\vec{d}}(w) \circ \psi(\alpha^{-1}) \circ \rho_{\vec{d}}(w_{\lambda}))(v) &= \sum_{\vec{0}\leq \vec{i}\leq \vec{d}} c_{\vec{i}} \sum_{\vec{0}\leq \vec{j}\leq \vec{i}} \binom{\vec{i}}{\vec{j}} \uni^{\vec{j}}(-\lambda)^{\vec{i}-\vec{j}} e_{\vec{d},\vec{j}} \\
&= \sum_{\vec{0}\leq \vec{j}\leq \vec{d}} \Biggl(\uni^{\vec{j}} \sum_{\vec{j}\leq \vec{i}\leq \vec{d}} c_{\vec{i}} \binom{\vec{i}}{\vec{j}} (-\lambda)^{\vec{i}-\vec{j}} \Biggr) e_{\vec{d},\vec{j}}.
\end{align*}
\end{proof}
\end{lemma}

\section{Un critère de séparation}\label{separ}



\subsection{Le résultat principal}\label{capo}

Conservons les notations des Sections \ref{lattici} et \ref{formulario} et fixons $\iota$ un plongement de $F$ dans $E$. Notons:
\[
S^+ = \{\sigma \in S, d_{\sigma} \neq 0 \}\subseteq S 
\]  
et, pour $l$ dans $\{0,\ldots, f-1 \}$ posons:
\[
J_l = \Big\{\sigma\in S^+, \sigma(\lambda) = \iota(\lambda)^{p^l} \ \forall \lambda \in I_1  \Big\},
\]
où rappelons que $I_1$ désigne l'ensemble des $[\zeta]$, $\zeta \in k_F$. En particulier, notons que
\[
\coprod_{l \in \{0,\ldots, f-1\}} J_l = S^+ \quad \mbox{et} \quad \forall\,  l \in \{0,\ldots, f-1\}, |J_l|\leq e.
\]  
Si $\sigma \in J_l$ on pose $\gamma_{\sigma} = l$. Si $i \in \mathbb{Z}$ on désigne par $\overline{i}$ l'unique élément de $\{1,\ldots,f \}$ tel que $i \equiv \overline{i}\, \mathrm{mod}\, f$ et posons pour $\sigma \in J_l$:
\[
v_{\sigma} = \inf \Big\{1\leq i \leq f, J_{\overline{l+i}} \neq \emptyset  \Big\}.
\]

Soit $a_p \in \mathfrak{p}_E$ et $\eta\colon \mathrm{Gal}(\overline{\mathbb{Q}}_p/F) \to \Oe^{\times}$ un caractère cristallin. On pose alors
\[
\Pi_{\vec{d},a_p,\eta} = \frac{\ind \rho_{\vec{d}}}{(T-a_p)(\ind \rho_{\vec{d}})}   \otimes (\eta \circ \mathrm{det}).
\]
Cette représentation est localement algébrique et peut se réaliser comme le produit tensoriel d'une représentation algébrique par une représentation lisse. Plus précisement on a le résultat suivant.

\begin{prop}\label{indotte}
Posons $u_{\sigma} = \frac{d_{\sigma}}{2}$ pour tout $\sigma \in S$. Avec les notations précédentes, on a:
\medbreak
\begin{itemize}
\item[(i)] Si $a_p \notin \{\pm(q \uni^{\vec{u}} + \uni^{\vec{u}}) \}$ alors  $\Pi_{\vec{d},a_p,\eta}$ est algébriquement irréductible et
\[
\Pi_{\vec{d},a_p,\eta} \simeq \rho_{\vec{d}} \otimes \mathrm{Ind}_P^G (\mathrm{nr}(\lambda_1^{-1})\otimes \mathrm{nr}(\lambda_2^{-1}))
\] 
où
\[
\lambda_1^f \lambda_2^f = q \uni^{\vec{d}}, \quad  \lambda_1^f + \lambda_2^f = a_p.
\]
\item[(ii)] Si $a_p \in \{\pm(q \uni^{\vec{u}} + \uni^{\vec{u}}) \}$ on a une suite exacte
\[
0 \to \rho_{\vec{d}} \otimes \mathrm{St}_G \otimes (\eta\, \mathrm{nr}(\delta)\circ \mathrm{det}) \to \Pi_{\vec{d},a_p,\eta} \to \rho_{\vec{d}} \otimes (\eta\, \mathrm{nr}(\delta) \circ \mathrm{det}) \to 0
\]
où  $\mathrm{St}_G = C^0(\mathbf{P}^1(F),E)/\{\mathrm{constantes} \}$ désigne la représentation de Steinberg de $G$, et $\delta = (q+1)/a_p$.
\end{itemize} 
\begin{proof}
Il s'agit d'une généralisation immédiate de \cite[Proposition 3.3]{breuilab}.
\end{proof}

\end{prop}

On suppose dorénavant $\eta = 1$ et on écrit $\Pi_{\vec{d},a_p}$ au lieu de $\Pi_{\vec{d},a_p,1}$. Comme dans \cite[§3.3]{breuilab} on définit:
\[
\Theta_{\vec{d},a_p} = \mathrm{Image} \Biggl(\ind \rho_{\vec{d}}^0 \longrightarrow \frac{\ind \rho_{\vec{d}}}{(T-a_p)(\ind \rho_{\vec{d}})} = \Pi_{\vec{d},a_p}  \Biggr),
\]
ou ce qui revient au même 
\begin{align}\label{analogia}
\Theta_{\vec{d},a_p} = \frac{\ind \rho_{\vec{d}}^0}{\ind \rho_{\vec{d}}^0\cap (T-a_p)(\ind \rho_{\vec{d}})}.
\end{align}
C'est un réseau au sens de la Définition \ref{lattice} et, puisque $\ind \rho_{\vec{d}}^0$ est de type fini en tant que $\Oe[G]$-module,   on déduit que $\Theta_{\vec{d},a_p}$ est aussi un $\Oe[G]$-module de type fini. De plus, comme 
\[
\forall h\in \ind \rho_{\vec{d}}^0, \quad T(h) \in \ind \rho_{\vec{d}}^0, 
\]
on déduit: 
\[
(T-a_p) ( \ind \rho_{\vec{d}}^0)  \subseteq     (T-a_p)(\ind \rho_{\vec{d}}) \cap \ind \rho_{\vec{d}}^0.
\]
Il découle alors immédiatement de \eqref{analogia} que l'on a une application surjective de $\Oe[G]$-modules:
\[
\theta\colon \frac{\ind \rho_{\vec{d}}^0}{(T-a_p)(\ind \rho_{\vec{d}}^0)} \onto \Theta_{\vec{d},a_p}.
\] 

Nous nous proposons de donner ici un critère pour que l'application $\theta$ soit injective. Commençons par deux lemmes techniques.

\begin{lemma}\label{tecniq}
Supposons qu'il existe $l\in \{0,\ldots, f-1\}$ tel que $|J_l|> 1$. Alors l'application $\theta$ n'est pas injective.
\begin{proof}
Notons que $\theta$ est injective si et seulement si l'on a l'inclusion suivante:
\begin{align}\label{inclusionelemma}
\big(T-a_p \big)\big(\ind \rho_{\vec{d}}\big) \cap \ind \rho_{\vec{d}}^0  \subseteq  \big(T-a_p \big) \big(\ind\rho_{\vec{d}}^0 \big).  
\end{align}
Il suffit alors de montrer que l'inclusion \ref{inclusionelemma} n'est pas vérifiée, où ce qui revient au même, que si $h$ est un élément dans $\ind \rho_{\vec{d}}$ tel que
\[
(T-a_p)(h) = T(h) - a_p h \in \ind \rho_{\vec{d}}^0,
\]
alors $h \notin \ind \rho_{\vec{d}}^0$. 
\medbreak
Par hypothèse il existe $l \in \{0,\ldots, f-1\}$, tel que $|J_l| > 1$. Autrement dit, il existe $\sigma, \tau \in S^+$ tels que $\sigma(\lambda) = \tau(\lambda)$ pour tout $\lambda \in I_1$. Notons $\vec{\alpha} = (\alpha_{\xi})_{\xi \in S}$ l'élément de $I_{\vec{d}}$ défini par:
\[
\alpha_{\xi} = \left\{ \begin{array}{ll}
1   & \mbox{si} \ \xi=\sigma, \\
0 & \mbox{si} \ \xi\neq \sigma
\end{array} \right.
\] 
et $v_1 = e_{\vec{d},\vec{\alpha}} \in \rho_{\vec{d}}$. Notons $\vec{\beta}= (\beta_{\xi})_{\xi \in S}$ l'élément de $I_{\vec{d}}$ défini par:
\[
\beta_{\xi} = \left\{ \begin{array}{ll}
1   & \mbox{si} \ \xi= \tau, \\
0 & \mbox{si} \ \xi \neq \tau
\end{array} \right.
\]
et $v_2= (-1) e_{\vec{d},\vec{\beta}} \in \rho_{\vec{d}}$. Posons $v = \delta^{-1} (v_1+v_2) \in \rho_{\vec{d}}$ où $\delta$ est celui des éléments $\iota(\uni)$, $a_p$ qui a la plus petite valuation et $h = [\mathrm{Id},v] \notin \ind \rho_{\vec{d}}^0$. Nous allons montrer  que $T(h)-a_p h \in \ind \rho_{\vec{d}}^0$. D'après le lemme \ref{tecnico1} on a:
\begin{align*}
T(h) - a_p h &= T^+(h)+T^-(h) - a_p h \\
&= \sum_{\lambda \in I_1} [g_{1,  \lambda }^0,(\rho_{\vec{d}}^0(w) \circ \psi(\alpha^{-1})\circ \rho_{\vec{d}}^0(w_{\lambda}))(v)] + [ \alpha ,\psi(\alpha^{-1})(v)]  - a_p [\mathrm{Id},v]
\end{align*}
et donc, comme $a_p [\mathrm{Id},v]\in \ind \rho_{\vec{d}}^0$,  il suffit de vérifier que les deux conditions suivantes sont satisfaites:
\begin{align}
[g_{1,  \lambda }^0,(\rho_{\vec{d}}^0(w) \circ \psi(\alpha^{-1})\circ \rho_{\vec{d}}^0(w_{\lambda}))(v)] &\in \ind \rho_{\vec{d}}^0 \quad \mbox{pour tout} \ \lambda \in I_1, \label{casus1} \\
[ \alpha ,\psi(\alpha^{-1})(v)] &\in \ind \rho_{\vec{d}}^0. \label{casus1.1}
\end{align}
Rappelons que pour tout $\xi \in S$, $U_{d_{\xi}}$ désigne la matrice diagonale définie dans \eqref{matricediag}. Posons $\varphi_{\lambda,\xi} = w \circ U_{d_{\xi}}\circ w_{\lambda}$ pour tout $\lambda \in I_1$. En utilisant la formule \eqref{azione3} on obtient:  
\begin{align*}
\varphi_{\lambda,\xi}(e_{d_{\xi},\alpha_{\xi}}) &= \left\{ \begin{array}{ll}
-\sigma(\lambda) e_{d_{\sigma},0}  + \sigma(\uni) e_{d_{\sigma},1}   & \quad \quad \ \mbox{si} \ \xi= \sigma, \\
e_{d_{\xi},0}  & \quad \quad \ \mbox{si} \ \xi\neq \sigma
\end{array} \right. \\
\varphi_{\lambda,\xi}(e_{d_{\xi},\beta_{\xi}}) &= \left\{ \begin{array}{ll}
 -\tau(\lambda) e_{d_{\tau},0}  + \tau(\uni) e_{d_{\tau},1}  & \,  \, \, \quad \quad \mbox{si} \ \xi= \tau, \\
e_{d_{\xi},0}    & \,  \, \, \quad \quad \mbox{si} \ \xi\neq \tau.
\end{array} \right.
\end{align*}
Posons $\varphi_{\lambda} = \rho_{\vec{d}}^0(w) \circ \psi(\alpha^{-1})\circ \rho_{\vec{d}}^0(w_{\lambda})$ pour tout $\lambda \in I_1$. Alors,   d'après la formule \eqref{azione4} on a:
\begin{align}
\forall \lambda \in I_1, \quad \varphi_{\lambda}(v_1) &\in -\sigma(\lambda) e_{\vec{d},\vec{0}}+ \sigma(\uni)\rho_{\vec{d}}^0 \label{casus7.3} \\ 
\forall \lambda \in I_1, \quad \varphi_{\lambda}(v_2) &\in \tau(\lambda)e_{\vec{d},\vec{0}}+\tau(\uni)\rho_{\vec{d}}^0 \label{casus7.4}
\end{align}
et comme pour tout $\lambda \in I_1$ on a $\sigma(\lambda) = \tau(\lambda)$ on déduit  de \eqref{casus7.3} et de \eqref{casus7.4}
\[
\forall \lambda \in I_1, \quad \varphi_{\lambda}(v) = \delta^{-1} \varphi_{\lambda}(v_1+v_2) \in (\delta^{-1}(\sigma(\uni)+\tau(\uni)))\rho_{\vec{d}}^0,
\]
d'où la condition \eqref{casus1}. Un calcul immédiat donne:
\begin{align}
U_{d_{\xi}}(e_{d_{\xi},\alpha_{\xi}}) &= \left\{ \begin{array}{ll}
\sigma(\uni)^{d_{\sigma}-1} e_{d_{\sigma},1}   &\quad \quad \mbox{si} \ \xi = \sigma, \\
\xi(\uni)^{d_{\xi}}e_{d_{\xi},0} &\quad \quad \mbox{si} \ \xi \neq \sigma
\end{array} \right. \label{tecs7.7} \\
U_{d_{\xi}}(e_{d_{\xi},\beta_{\xi}}) &= \left\{ \begin{array}{ll}
\tau(\uni)^{d_{\tau}-1} e_{d_{\tau},1}    &\quad\quad \,   \mbox{si} \ \xi= \tau, \\
\xi(\uni)^{d_{\xi}}e_{d_{\xi},0} &\quad\quad  \, \mbox{si} \ \xi\neq  \tau.
\end{array} \right. \label{tecs7}
\end{align}
Comme $|J_l|> 1$ on déduit de \eqref{tecs7.7} (resp. \eqref{tecs7}) que $\psi(\alpha^{-1})(v_1) \in \tau(\uni)\rho_{\vec{d}}^0$ (resp. $\psi(\alpha^{-1})(v_2) \in \sigma(\uni)\rho_{\vec{d}}^0$). On obtient 
\[
\psi(\alpha^{-1})(v) = \delta^{-1} \psi(\alpha^{-1})(v_1+v_2) \in (\delta^{-1}(\sigma(\uni)+ \tau(\uni)))\rho_{\vec{d}}^0 \subseteq \rho_{\vec{d}}^0,
\] 
d'où la condition \eqref{casus1.1}.

\end{proof}
\end{lemma}

\begin{lemma}\label{tecniq2}
Supposons qu'il existe $\sigma \in J_l$ tel que
\[
d_{\sigma}+1 \leq p^{v_{\sigma}}.
\] 
Alors l'application $\theta$ n'est pas injective.
\begin{proof}
Le même raisonnement que dans la démonstration du Lemme \ref{tecniq} montre qu'il suffit de construire explicitement un élément $h \in \ind\rho_{\vec{d}}$ tel que 
\[
(T-a_p)(h) \in \ind \rho_{\vec{d}}^0 \quad \mbox{et} \quad h\notin \ind\rho_{\vec{d}}^0.
\] 
On distingue trois cas.
\begin{itemize}
\item[(1)] $|S^+|\geq 2$ et $d_{\sigma} \geq p^{v_{\sigma}}$; 
\item[(2)] $|S^+| = 1$ et $d_{\sigma} \geq p^{v_{\sigma}} + 1$;
\item[(3)] $|S^+| = 1$ et $d_{\sigma} = p^{v_{\sigma}}$;
\end{itemize}
\medbreak
Cas (1). On peut supposer $v_{\sigma} < f$ car sinon, par hypothèse, il existe $l \in \{0,\ldots, f-1\}$ tel que $|J_l|\geq 2$. Or ce cas a déjà été  traité dans le Lemme \ref{tecniq}. Il existe alors $\tau \in S^+$ tel que 
\[
\forall \lambda \in I_1, \quad \sigma(\lambda)^{p^{v_{\sigma}}} = \tau(\lambda).
\]  
Notons $\vec{\alpha} = (\alpha_{\xi})_{\xi \in S}$ l'élément de $I_{\vec{d}}$ défini par:
\[
\alpha_{\xi} = \left\{ \begin{array}{ll}
p^{v_{\sigma}}   & \mbox{si} \ \xi = \sigma, \\
0 & \mbox{si} \ \xi \neq \sigma
\end{array} \right.
\] 
et $v_1 = (-1)^{p^{v_{\sigma}}} e_{\vec{d},\vec{\alpha}} \in \rho_{\vec{d}}$. Notons $\vec{\beta} = (\beta_{\xi})_{\xi \in S}$ l'élément de $I_{\vec{d}}$ défini par:
\[
\beta_{\xi} = \left\{ \begin{array}{ll}
1   & \mbox{si} \ \xi = \tau, \\
0 & \mbox{si} \ \xi \neq \tau
\end{array} \right.
\]
et $v_2= e_{\vec{d},\vec{\beta}} \in \rho_{\vec{d}}$. Posons $v = \delta^{-1} (v_1+v_2) \in \rho_{\vec{d}}$ où $\delta$ est celui des éléments $\iota(\uni)$, $a_p$ qui a la plus petite valuation et $f = [\mathrm{Id},v] \notin \ind \rho_{\vec{d}}^0$. Montrons  que $T(f)-a_pf \in \ind \rho_{\vec{d}}^0$. Comme déjà remarqué dans la preuve du Lemme \ref{tecniq} on est ramené à vérifier les deux conditions suivantes:    
\begin{align}
[g_{1,  \lambda }^0,(\rho_{\vec{d}}^0(w) \circ \psi(\alpha^{-1})\circ \rho_{\vec{d}}^0(w_{\lambda}))(v)] &\in \ind \rho_{\vec{d}}^0 \quad \mbox{pour tout} \ \lambda \in I_1, \label{casus2} \\
[ \alpha ,\psi(\alpha^{-1})(v)] &\in \ind \rho_{\vec{d}}^0. \label{casus2.1}
\end{align}
Posons $\varphi_{\lambda,\xi} = w \circ U_{d_{\xi}}\circ w_{\lambda}$ pour tout $\lambda \in I_1$ et tout $\xi \in S$. D'après la formule \eqref{azione3} on a:  
\begin{align*}
\varphi_{\lambda,\xi}(e_{d_{\xi},\alpha_{\xi}}) &= \left\{ \begin{array}{ll}
(-1)^{p^{v_{\sigma}}} \sigma(\lambda)^{p^{v_{\sigma}}} e_{d_{\sigma},0}+ \sigma(\uni) \sum_{l=1}^{p^{v_{\sigma}}}\binom{p^{v_{\sigma}}}{l}\sigma(\uni)^{l-1} \sigma(-\lambda)^{p^{v_{\sigma}}-l} e_{d_{\sigma},l}    & \mbox{si} \ \xi = \sigma, \\
e_{d_{\xi},0} & \mbox{si} \ \xi \neq \sigma
\end{array} \right. \\
\varphi_{\lambda,\xi}(e_{d_{\xi},\beta_{\xi}}) &= \left\{ \begin{array}{ll}
-\tau(\lambda)e_{d_{\tau},0}+\tau(\uni) e_{d_{\tau},1}    & \quad \quad \quad\quad \quad \quad\quad \; \, \ \quad \quad \quad \quad \quad \quad \quad \quad \quad \quad \mbox{si} \ \xi= \tau, \\
e_{d_{\xi},0} & \quad \quad\quad \quad \quad \quad \quad \quad \quad \quad \quad \quad  \quad \ \; \, \quad \quad \quad \quad \mbox{si} \ \xi \neq  \tau.
\end{array} \right.
\end{align*}
Posons $\varphi_{\lambda} = \rho_{\vec{d}}^0(w) \circ \psi(\alpha^{-1})\circ \rho_{\vec{d}}^0(w_{\lambda})$ pour tout $\lambda \in I_1$. Alors,  d'après la formule \eqref{azione4} on a:
\begin{align}
\forall \lambda \in I_1, \quad \varphi_{\lambda}(v_1) &\in \sigma(\lambda)^{p^{v_{\sigma}}} e_{\vec{d},\vec{0}}+ \sigma(\uni)\rho_{\vec{d}}^0 \label{casus2.3} \\ 
\forall \lambda \in I_1, \quad \varphi_{\lambda}(v_2) &\in -\tau(\lambda)e_{\vec{d},\vec{0}}+ \tau(\uni)\rho_{\vec{d}}^0 \label{casus2.4}
\end{align}
et comme pour tout $\lambda \in I_1$ on a $\sigma(\lambda)^{p^{v_{\sigma}}} = \tau(\lambda)$ on déduit  de \eqref{casus2.3} et de \eqref{casus2.4}
\[
\forall \lambda \in I_1, \quad \varphi_{\lambda}(v) = \delta^{-1} \varphi_{\lambda}(v_1+v_2) \in (\delta^{-1}\iota(\uni))\rho_{\vec{d}}^0,
\]
d'où la condition \eqref{casus2}. Un calcul immédiat donne:
\begin{align}
U_{d_{\xi}}(e_{d_{\xi},\alpha_{\xi}}) &= \left\{ \begin{array}{ll}
\sigma(\uni)^{d_{\sigma}-p^{v_{\sigma}}} e_{d_{\sigma},p^{v_{\sigma}}}   & \mbox{si} \ \xi = \sigma, \\
\xi(\uni)^{d_{\xi}}e_{d_{\xi},0} & \mbox{si} \ \xi \neq \sigma
\end{array} \right. \label{tecs1} \\
U_{d_{\xi}}(e_{d_{\xi},\beta_{\xi}}) &= \left\{ \begin{array}{ll}
\tau(\uni)^{d_{\tau}-1}e_{d_{\tau},1}    & \quad \ \, \  \mbox{si} \ \xi= \tau, \\
\xi(\uni)^{d_{\xi}} e_{d_{\xi},0} & \quad\ \, \ \mbox{si} \ \xi \neq \tau.
\end{array} \right. \label{tecs2}
\end{align}
Comme $|S^+|> 1$ (et donc $v_{\sigma} < f$) on déduit de \eqref{tecs1} que $\psi(\alpha^{-1})(v_1) \in \iota(\uni)\rho_{\vec{d}}^0$. Comme $d_{\sigma}\geq p^{v_{\sigma}}$ on déduit de \eqref{tecs2} que $\psi(\alpha^{-1})(v_2) \in \iota(\uni)\rho_{\vec{d}}^0$. D'après ce qui précède on a: 
\[
\psi(\alpha^{-1})(v) = \delta^{-1} \psi(\alpha^{-1})(v_1+v_2) \in (\delta^{-1}\iota(\uni))\rho_{\vec{d}}^0 \subseteq \rho_{\vec{d}}^0,
\] 
d'où la condition \eqref{casus2.1}.
\medbreak
Cas (2). Puisque $|S^+| = 1$ on a $v_{\sigma} = f$ et donc
\[
d_{\sigma} \geq p^{v_{\sigma}}+1 = p^f + 1 = q+1.
\]
Notons que dans ce cas $\rho_{\vec{d}} = (\mathrm{Sym}^{d_{\sigma}}E^2)^{\sigma}$ et de même $\rho_{\vec{d}}^0 = (\mathrm{Sym}^{d_{\sigma}}\Oe^2)^{\sigma}$. Posons: 
\[
v = \delta^{-1} ( e_{d_{\sigma},1} + (-1)^{q} e_{d_{\sigma},q}) \in  \rho_{\vec{d}}
\] 
où $\delta$ est celui des éléments $\iota(\uni)$, $a_p$ qui a la plus petite valuation et $h = [\mathrm{Id},v] \notin \ind \rho_{\vec{d}}^0$. Un raisonnement analogue à celui donné dans le cas (1) (dont on laisse les détails au lecteur) montre que $T(h)-a_p h \in \ind \rho_{\vec{d}}^0$. 
\medbreak 
Cas (3). On a supposé 
\[
|S^+|=1, \quad  d_{\sigma} = p^{v_{\sigma}} = p^f = q.
\]
Posons:
\begin{align*}
\forall \lambda \in I_1, \quad v_{2,\lambda} &=  \left\{ \begin{array}{ll}
\delta^{-1}\sigma(\lambda)^{q-2}(-e_{q,1}-e_{q,q})   & \; \mbox{si} \  p=2, \\
\delta^{-1}\sigma(\lambda)^{q-2}(e_{q,1}-e_{q,q}) & \; \mbox{si}  \ p \neq 2,
\end{array} \right.    \\
                         v_{0}   &= \left\{ \begin{array}{ll}
\delta^{-1}(q-1)(e_{q,0}+e_{q,q-1})   & \,\; \mbox{si} \ p=2, \\
\delta^{-1}(q-1)(e_{q,0}-e_{q,q-1}) & \, \;  \mbox{si} \ p \neq 2,
\end{array} \right. \\
                         h_2 &= \sum_{\lambda \in I_1} [g_{2,\uni\lambda}^0,v_{2,\lambda}], \\
                         h_0 &= [\mathrm{Id},v_{0}],
\end{align*}
où $\delta$ est celui des éléments $\sigma(\uni)$, $a_p$ qui a la plus petite valuation et notons $h = h_0 + h_2 \notin \ind \rho_{\vec{d}}^0$. Plaçons-nous dans le cas où $q$ est la puissance d'un nombre premier impair (l'autre cas se traitant de manière analogue). Nous allons montrer que $T(h)-a_p h \in \ind \rho_{\vec{d}}^0$. D'après le Lemme \ref{tecnico1} la fonction $T(h)-a_p h$ peut s'écrire somme des quatre fonctions qui ont des supports deux à deux disjoints:
\[
T(h) - a_p h = T^- (h_0) + (T^+(h_0) + T^-(h_2))+ T^+(h_2) - a_p h
\]
et donc, pour montrer que $T(h) - a_p h \in \ind \rho_{\vec{d}}^0$ il faut voir que chaque fonction est dans $\ind \rho_{\vec{d}}^0$. D'après le Lemme \ref{tecnico1} on a:
\begin{align*}
T^-(h_0) &= [\alpha,\psi(\alpha^{-1})(v_0)] \\
&= [\alpha,\delta^{-1}(q-1)\sigma(\uni) (\sigma(\uni)^{q-1}e_{q,0}-e_{q,q-1})] \in (\delta^{-1}\sigma(\uni))\ind \rho_{\vec{d}}^0 \subseteq \ind \rho_{\vec{d}}^0.
\end{align*}
Posons: 
\begin{align*}
\forall \mu \in I_1, \quad \varphi_{\mu}&=\rho_{\vec{d}}^0(w) \circ \psi(\alpha^{-1})\circ \rho_{\vec{d}}^0(w_{\mu}), \\ 
\forall \mu \in I_1, \quad \phi_{\mu}&=\rho_{\vec{d}}^0(w) \circ \rho_{\vec{d}}^0(w_{-\mu}) \circ \psi(\alpha^{-1}).
\end{align*}
Par linéarité de l'opérateur $T^+$ et d'après le Lemme \ref{tecnico1} on a:
\begin{align*}
T^+(h_2) \in \ind \rho_{\vec{d}}^0 &\Leftrightarrow \ \forall \lambda \in I_1, \ T^+([g_{2,\uni \lambda}^0,v_{2,\lambda}]) \in \ind \rho_{\vec{d}}^0 \\
&\Leftrightarrow \ \forall \lambda \in I_1, \  \forall \mu \in I_1, \ \varphi_{\mu}(v_{2,\lambda}) \in  \rho_{\vec{d}}^0.
\end{align*}
D'après la formule \eqref{azione3} on a pour tout $\mu \in I_1$ et tout $\lambda \in I_1$:
\begin{align*}
\varphi_{\mu}(v_{2,\lambda}) &=  \varphi_{\mu}(\delta^{-1}\sigma(\lambda)^{q-2}(e_{q,1}-e_{q,q})) \\
 &= \delta^{-1}\sigma(\lambda)^{q-2} \Big( \sigma(\uni) e_{q,1} - \sigma(\uni) \sum_{\alpha = 1}^q \binom{q}{\alpha}(-1)^{q-\alpha} \sigma(\uni)^{\alpha - 1} \sigma(\mu)^{q-\alpha} e_{q,\alpha}\Big) \in \rho_{\vec{d}}^0,
\end{align*} 
(pour ce calcul on a utilisé le fait que $q$ est une puissance d'un premier impair). Il nous reste à montrer que $T^+(h_0) + T^-(h_2) \in \ind \rho_{\vec{d}}^0$. Remarquons que d'après le Lemme \ref{tecnico1} on a:
\begin{align*}
T^-(h_2) &= \sum_{\lambda \in I_1} T^-([g_{2,\uni\lambda}^0,v_{2,\lambda}]) = \sum_{\lambda \in I_1} [g_{1,0}^0,\phi_{\lambda}(v_{2,\lambda})], \\
T^+(h_0) &=   T^+([\mathrm{Id},v_0]) = \sum_{\lambda \in I_1}[g_{1,\lambda}^0,\varphi_{\lambda}(v_0)]. 
\end{align*}
En utilisant la formule \eqref{azione3} on a:
\[
\forall \lambda \in I_{1}\backslash \{0\}, \quad \varphi_{\lambda}(v_0) = \delta^{-1}(q-1)\sigma(\uni) \Big(\sum_{\alpha = 1}^{q-1} \binom{q-1}{\alpha}\sigma(\uni)^{\alpha-1} \sigma(-\lambda)^{q-1-\alpha}e_{q,\alpha}  \Big) \in \rho_{\vec{d}}^0
\]
(pour ce calcul on a de même utilisé le fait que $q$ est une puissance d'un premier impair) ce qui implique:
\[
\forall \lambda \in I_{1}\backslash \{0\}, \quad [g_{1,\lambda}^0,\varphi_{\lambda}(v_0)] \in \ind \rho_{\vec{d}}^0.
\]
Il nous reste alors à calculer
\begin{align*}
\varphi_{0}(v_0) + \sum_{\lambda \in I_1} \phi_{\lambda}(v_{2,\lambda}).
\end{align*}
En utilisant la formule \eqref{azione3} on obtient:
\begin{align*}
\phi_{\lambda}(v_{2,\lambda}) &\in \delta^{-1} (-\sigma(\lambda)^{2q-2}e_{q,0}-\sigma(\lambda)^{q-2}e_{q,q}) + (\delta^{-1}\sigma(\uni))\rho_{\vec{d}}^0\\
\varphi_0(v_0) &\in \delta^{-1} ((q-1)e_{q,0}) + (\delta^{-1}\sigma(\uni))\rho_{\vec{d}}^0,
\end{align*}
et donc, d'après les relations 
\[
\sum_{\lambda \in I_1}\sigma(\lambda)^{2q-2} = q-1, \quad p| \sum_{\lambda \in I_1}\sigma(\lambda)^{q-2} 
\]
on déduit: 
\[
\varphi_{0}(v_0) + \sum_{\lambda \in I_1} \phi_{\lambda}(v_{2,\lambda}) \in \rho_{\vec{d}}^0, 
\] 
ce qui permet de conclure.

\end{proof}
\end{lemma}

Le théorème suivant fournit deux conditions nécessaires et suffisantes sur le vecteur $\vec{d} = (d_{\sigma})_{\sigma \in S}$ pour que l'application $\theta$ soit injective.

\begin{theo}\label{principe}
Avec les notations précédentes l'application $\theta$ est injective (et donc un isomorphisme) si et seulement si les deux conditions suivantes sont satisfaites: 
\begin{itemize}
\item[(i)] Pour tout $l\in \{0,\ldots, f-1\}$ on a $|J_l|\leq 1$;
\item[(ii)] Si $\sigma \in J_l$ on a  
\[
d_{\sigma}+1 \leq p^{v_{\sigma}}.
\]
\end{itemize}
\begin{proof}
Comme déjà remarqué dans la preuve du Lemme \ref{tecniq}, l'application $\theta$ est injective si et seulement si l'on a l'inclusion suivante:
\begin{align}\label{inclusione}
\big(T-a_p\big)\big(\ind \rho_{\vec{d}}\big) \cap \ind \rho_{\vec{d}}^0  \subseteq  \big(T-a_p\big) \big(\ind\rho_{\vec{d}}^0 \big).  
\end{align}

La preuve comporte deux étapes.
\medbreak
\begin{itemize}
\item[(1)] On suppose que les conditions (i) et (ii) soient satisfaites. On montre que l'inclusion \eqref{inclusione} est vérifiée, où ce qui revient au même, que si $h$ est un élément dans $\ind \rho_{\vec{d}}$ tel que
\[
(T-a_p)(h) = T(h) - a_p h \in \ind \rho_{\vec{d}}^0,
\]
alors $h$ est dans $\ind \rho_{\vec{d}}^0$.
\item[(2)] On suppose que (i) ou bien (ii) ne soit pas satisfaite. On construit explicitement un $h$ dans $\ind \rho_{\vec{d}}$ tel que
\[
(T-a_p)(h) = T(h) - a_p h \in \ind \rho_{\vec{d}}^0 \quad \mbox{et} \quad h\notin \ind \rho_{\vec{d}}^0.
\]
\end{itemize}
\bigbreak
(1). Supposons que les conditions (i) et (ii) soient satisfaites. Soit $h \in \ind \rho_{\vec{d}}$ tel que:
\begin{align} \label{condizione}
(T-a_p)(h) = T(h) - a_p h \in \ind \rho_{\vec{d}}^0,
\end{align}
et notons $n$ le plus petit entier tel que $h \in B_n$ et écrivons $h = \sum_{m=0}^n h_m$ où $h_m \in S_m$. On en déduit:
\begin{align*}
T(h)-a_p h &= T^+(h) + T^-(h)- a_p h \\
&= T^+(h_n) + T^+(h-h_n)+ T^-(h)- a_p h, 
\end{align*}
avec $T^+(h_n) \in S_{n+1}$ et $T^+(h-h_n)+ T^-(h)- a_p h \in B_n$. L'hypothèse \eqref{condizione} implique  que $T^+(h_n) \in \ind \rho_{\vec{d}}^0$. Montrons que cela implique que $h_n \in \ind \rho_{\vec{d}}^0$. Il suffit de le montrer pour $h_n$ de la forme $[g_{n,\mu}^0,v]$, car, en appliquant $\beta$, on le déduit pour $h_n$ de la forme $[g_{n,\mu}^1,w]$, et enfin, par linéarité et en vertu du Corollaire \ref{tecnico4}, pour n'importe quel $h_n$. D'après le Lemme \ref{tecnico1} on a:
\begin{align}\label{tipi}
T^+([g_{n,\mu}^0,v]) \in \ind \rho_{\vec{d}}^0 \Leftrightarrow  (\rho_{\vec{d}}^0(w) \circ \psi(\alpha^{-1})\circ \rho_{\vec{d}}^0(w_{\lambda}))(v) \in V_{\rho_{\vec{d}}^0} \quad \mbox{pour tout} \ \lambda \in I_1.
\end{align}
Notons $v = \sum_{\vec{0}\leq \vec{i}\leq \vec{d}} c_{\vec{i}} e_{\vec{d},\vec{i}}$ avec $c_{\vec{i}} \in E$. D'après le Lemme \ref{tecnico5} on a:
\begin{align*}
(\rho_{\vec{d}}^0(w) \circ \psi(\alpha^{-1})\circ \rho_{\vec{d}}^0(w_{\lambda}))(v) \in V_{\rho_{\vec{d}}^0} \Leftrightarrow \uni^{\vec{j}} \sum_{\vec{j}\leq \vec{i}\leq \vec{d}} c_{\vec{i}} \binom{\vec{i}}{\vec{j}} (-\lambda)^{\vec{i}-\vec{j}} \in \Oe, \quad \vec{0}\leq \vec{j}\leq \vec{d}, \ \lambda \in I_1
\end{align*}
et donc, en particulier, pour $\vec{j}=\vec{0}$ on déduit de \eqref{tipi} l'implication suivante: 
\begin{align}\label{implic}
\forall \lambda \in I_1,\quad T^+([g_{n,\mu}^0,v]) \in \ind \rho_{\vec{d}}^0 \Rightarrow \sum_{\vec{0}\leq \vec{i}\leq \vec{d}} c_{\vec{i}}  (-\lambda)^{\vec{i}} \in \Oe.
\end{align}
Remarquons que pour $\lambda = 0$ on déduit que $c_{\vec{0}} \in \Oe$ et donc l'implication \eqref{implic} est équivalente à l'implication suivante:
\[
\forall \lambda \in k_F^{\times}, \quad T^+([g_{n,\mu}^0,v]) \in \ind \rho_{\vec{d}}^0 \Rightarrow \sum_{\substack{\vec{0}\leqslant \vec{i}\leqslant \vec{d}\\ \vec{i}\neq \vec{0}}} c_{\vec{i}}  (-[\lambda])^{\vec{i}} \in \Oe.
\] 
Soit $\zeta$ un générateur du groupe cyclique $k_F^{\times}$ et posons pour tout $\sigma \in S$, $u_{\sigma} = d_{\sigma}+1$. Alors on a:
\[
\forall \lambda \in k_F^{\times}, \quad \sum_{\substack{\vec{0}\leqslant \vec{i}\leqslant \vec{d}\\ \vec{i}\neq \vec{0}}} c_{\vec{i}}  (-[\lambda])^{\vec{i}} \in \Oe \Rightarrow   \sum_{\substack{\vec{0}\leqslant \vec{i}\leqslant \vec{d}\\ \vec{i}\neq \vec{0}}} c_{\vec{i}}  ([\zeta^j])^{\vec{i}} \in \Oe, \ \ 0\leq j \leq \Big(\prod_{\sigma \in S} u_{\sigma}\Big)-1.
\]
On obtient un système de $\prod_{\sigma \in S} u_{\sigma}$ équations linéaires à $\prod_{\sigma \in S} u_{\sigma}$  inconnues dont la matrice du système homogène associée est donnée par:
\[
A= \begin{bmatrix}
1 \\
[\zeta]^{\vec{i}} \\
[\zeta^2]^{\vec{i}} \\
\vdots &   \\
[\zeta^{\prod_{\sigma \in S} u_{\sigma} -1}]^{\vec{i}} &  
\end{bmatrix}_{\substack{\vec{0}\leqslant \vec{i}\leqslant \vec{d}\\ \vec{i}\neq \vec{0}}}
\]
Puisque $A$ est une matrice de Vandermonde, cela implique que son déterminant est égal à:
\begin{align}\label{vand}
\prod_{\vec{i}\prec \vec{j}, \vec{i}\neq \vec{j}}  ([\zeta]^{\vec{j}}- [\zeta]^{\vec{i}}) 
\end{align}
(rappelons que l'on a muni $I_{\vec{d}}$ de l'ordre lexicographique $\prec$). Comme pour tout $\sigma \in S^+$ on a $\sigma([\zeta]) = \iota([\zeta])^{p^{\gamma_{\sigma}}}$ où $0\leq \gamma_{\sigma} \leq f-1$, alors pour tout $\vec{j} \in I_{\vec{d}}$ on a par définition:
\begin{align*}
[\zeta]^{\vec{j}} = \prod_{\sigma \in S^+} \sigma([\zeta])^{j_{\sigma}} = \iota([\zeta])^{\sum_{\sigma \in S^+} j_{\sigma}p^{\gamma_{\sigma}}}
\end{align*}
et donc, d'après les hypothèses (i) et (ii), on obtient de manière naturelle une application injective:
\[
\big\{[\zeta]^{\vec{j}}, \vec{0}\leq \vec{j} \leq \vec{d} \ \mbox{et} \ \vec{j} \neq \vec{0}   \big\} \into \big\{\iota([\zeta])^{\alpha}, 0\leq \alpha \leq p^f-2   \big\} 
\]
qui est bijective si $d_{\sigma} +1 = p^{v_{\sigma}}$ pour tout $\sigma \in S^+$. Or, si $0\leq \alpha < \beta \leq q-2$ alors $\iota([\zeta])^{\alpha}- \iota([\zeta])^{\beta} \in \Oe^{\times}$ et donc, en utilisant \eqref{vand}, on déduit que $\mathrm{det}(A) \in \Oe^{\times}$. Puisque par ce qui précède $c_{\vec{i}} \in \Oe$ pour tout $\vec{i} \in I_{\vec{d}}$, ou ce qui revient au même $v \in \rho_{\vec{d}}^0$, on déduit que $h_n\in \ind \rho_{\vec{d}}^0$. En remplaçant $h$ par $h-h_n$, le même raisonnement montre que $h_{n-1} \in \ind \rho_{\vec{d}}^0$. Une récurrence évidente donne $h_i \in \ind \rho_{\vec{d}}^0$ pour tout $i \in \{0,\ldots, n\}$ et donc $h \in \ind \rho_{\vec{d}}^0$, d'où la première étape. 
\medbreak
(2). C'est une conséquence immédiate des Lemmes \ref{tecniq} et \ref{tecniq2}.

\end{proof}
\end{theo}

Soit $a_p \in \Oe$ et supposons que le vecteur d'entiers $\vec{d}$ vérifie les conditions (i) et (ii)  du Théorème \ref{principe}, c'est-à-dire:
\begin{itemize}
\item[(i)] Pour $l \in \{0,\ldots, f-1\}$ on a $|J_l|\leq 1$;
\item[(ii)] Si $\sigma \in J_l$ on a  
\[
d_{\sigma}+1 \leq p^{v_{\sigma}}.
\]
\end{itemize}  
Notons pour tout $n \in \mathbb{Z}_{\geq 0}$:
\[
B_n(E) = \{g \in \ind \rho_{\vec{d}}, g\in B_{n} \}.
\]
Une conséquence simple mais intéressante du Théorème \ref{principe} est la proposition suivante, démontrée initialement par Breuil pour $F=\Q$ et $d\leq 2p$ si $p\neq 2$ et $k<4$ si $p=2$ dans \cite[Théorème 4.1 et Corollaire 4.2]{breuilab}. 

\begin{prop}\label{libero}
Le $\Oe$-réseau $\Theta_{\vec{d},a_p}$ est séparé au sens de la Définition \ref{lattice}.	
\begin{proof}
L'argument de \cite[Corollaire 4.2]{breuilab} s'étend sans problème. Il suffit de montrer que $\Theta_{\vec{d},a_p}$ ne contient pas de $E$-droite. Cela revient à montrer que si $h \in \ind \rho_{\vec{d}}^0$ est tel qu'il existe des $h_n \in \ind \rho_{\vec{d}}$ vérifiant:
\begin{align}\label{latte}
\forall n \in \mathbb{Z}_{\geq 0},\quad h - (T-a_p)(h_n) \in p^n \big(\ind \rho_{\vec{d}}^0\big) 
\end{align}
alors $h \in (T-a_p)(\ind \rho_{\vec{d}})$. Notons $N$ (resp. $N'$) le plus petit entier positif ou nul tel que $h \in B_N$ (resp. $h_n \in B_{N'}$) et écrivons $h_{n} = \sum_{m=0}^{N'} h_{n,m}$ où $h_{n,m} \in S_m$. Si $N'\geq N$ la relation \eqref{latte} implique que $T^+(h_{n,N'}) \in p^n (\ind \rho_{\vec{d}}^0)$ et donc, la même preuve que celle donnée dans le Théorème \ref{principe} ($\Rightarrow$) implique que $h_{n,N'} \in p^n (\ind \rho_{\vec{d}}^0)$. Une récurrence descendante immédiate montre alors que $h_{n,m} \in p^n (\ind \rho_{\vec{d}}^0)$ pour tout $N\leq m\leq N'$. En résumant on a:
\[
h_n \in B_{N-1}(E) + p^n \big(\ind \rho_{\vec{d}}^0\big). 
\]
Ainsi:
\[
\forall n \in \mathbb{Z}_{\geq 0}, \quad h \in (T-a_p) (B_{N-1}(E)) + p^n \big(\ind \rho_{\vec{d}}^0\big) 
\]
et, comme $(T-a_p) ( B_{N-1}(E))$ est un $E$-espace vectoriel complet pour la topologie $p$-adique car de dimension finie, on en déduit $h \in (T-a_p) ( B_{N-1}(E))$.

\end{proof}
\end{prop}

\subsection{Conséquences}\label{conseguenz}

Conservons les notations de la Section \ref{capo}. On fixe
\smallbreak
\begin{itemize}
\item[$\bullet$] $(\alpha,\beta) \in E^{\times} \times E^{\times}$.
\item[$\bullet$] un $|S|$-uplet d'entiers positifs $\vec{d}$. 
\end{itemize}
Rappelons que $\rho_{\vec{d}}$ désigne la représentation algébrique irréductible de $\GL(F)$ introduite dans le Paragraphe \ref{repr}. 
\medbreak
Nous allons ici rappeler, dans un cadre particulier, une conjecture formulée par Breuil et Schneider (\cite{bs}). Nous montrons ensuite comment utiliser le résultat principal de la Section \ref{capo} (Théorème \ref{principe}) pour donner, dans quelques cas, une réponse positive à cette conjecture.

\medbreak
Notons $(r,V)$ la $E$-représentation du groupe de Weil-Deligne de $F$ de dimension $2$, reductible, non ramifiée et qui envoie le Frobenius arithmétique sur la matrice
\[
\begin{bmatrix} {\alpha} & {0} \cr {0} & {\beta} \end{bmatrix}.
\]
Notons $\pi^{unit}$ la représentation lisse et irréductible de $G$ sur $E$ qui correspond à $(r,V)$ via la correspondance de Langlands locale, normalisée de sorte que son caractère central soit $\mathrm{det}(r,V)\circ \mathrm{rec}^{-1}$. Dans \cite[§4]{bs} est décrite une construction pour associer à $\pi^{unit}$ une représentation lisse de $G$ sur $E$ que l'on note $\pi$. Cette construction tient compte du fait que $\pi^{unit}$ soit générique ou non. Dans ce cadre particulier on obtient dans les deux cas:
\[
\pi = \Ind_P^G (\mathrm{nr}(\alpha^{-1})\otimes \mathrm{nr}(p\beta^{-1})).
\]  

\begin{rem}{\rm
Par définition les caractères non ramifiés $\mathrm{nr}(\alpha^{-1})$ et $\mathrm{nr}(p\beta^{-1})$ dépendent seulement de $\alpha^f$ et de $\beta^f$. Notons que $\pi$ est une induite parabolique lisse et non ramifiée qui est irréductible si $(\alpha\beta^{-1})^f \neq q$ et $(\alpha\beta^{-1})^f \neq q^{-1}$. Si $(\alpha\beta^{-1})^f = q$ (resp. $(\alpha\beta^{-1})^f \neq q^{-1}$) alors $\pi$ est la torsion par $\mathrm{nr}(\beta^{-1})\circ \mathrm{det}$ (resp. $\mathrm{nr}(\alpha^{-1})\circ \mathrm{det}$) de l'unique extension non scindée de la représentation triviale par la Steinberg (resp. de la Steinberg par la représentation triviale). }
\end{rem}

Appelons $\varphi$-module la donnée d'un $F_0 \otimes_{\Q}E$-module libre de rang fini $D$ muni d'un automorphisme $F_0$-semi-linéaire (par rapport au Frobenius sur $F_0$) et $E$-linéaire $\varphi$. Notons que $\varphi^f$ est une application $F_0\otimes_{\Q}E$-linéaire 
 et que l'isomorphisme $F_0\otimes_{\Q}E \simeq \prod_{\sigma_0\colon F_0\into E} E$, $h\otimes e \mapsto (\sigma_0(h)e)_{\sigma_0}$ induit un isomorphisme
\[
D \simeq \prod_{\sigma_0\colon F_0 \into E} D_{\sigma_0}
\] 
où $D_{\sigma_0} = (0,0,\ldots,0,1_{\sigma_0},0,\ldots,0)\cdot D$. On définit 
\[
t_N(D) = \frac{1}{[F:\Q]} val_F(\mathrm{det}_{F_0}(\varphi^f|_D)).
\]

Si $D$ est un $\varphi$-module on peut lui associer de manière explicite une représentation de Weil-Deligne par la méthode décrite dans \cite{fon}. Plus précisement, choisissons un plongement $\sigma_0\colon F_0 \into E$ et posons $U = D_{\sigma_0}$. Si $w\in W(\overline{\mathbb{Q}}_p/F)$, on définit $s(w) = \varphi^{-\alpha(w)}$ où $\alpha(w) \in f\mathbb{Z}$ désigne l'unique entier tel que l'action induite de $w$ sur $\overline{\mathbb{F}}_p$ soit la $\alpha(w)$-puissance du Frobenius arithmétique $x\mapsto x^p$. On vérifie que $s(w)\colon D \to D$ est $F_0 \otimes_{\Q} E$-linéaire, et donc induit un morphisme $E$-linéaire $s(w)\colon U \to U$. Le couple $(s,U)$ est une $E$-représentation non ramifiée du groupe de Weil-Deligne de $F$ qui ne dépend pas du choix de $\sigma_0$, mais à isomorphisme non canonique près (\cite[Lemme 2.2.1.2]{bm}). On le note $\mathrm{WD}(\varphi,D)$ et on note $\mathrm{WD}(\varphi,D)^{ss}$ sa $F$-semisimplification (\cite[§8.5]{deli}).
\medbreak
Soit $D$ un $\varphi$-module. Si on pose $D_F = D \otimes_{F_0}F$ alors l'isomorphisme $F\otimes_{\Q}E \simeq \prod_{\sigma\colon F\into E} E$, $h\otimes e \mapsto (\sigma(h)e)_{\sigma_0}$ induit un isomorphisme
\[
D \simeq \prod_{\sigma\colon F \into E} D_{\sigma}
\] 
où $D_{F,\sigma} = (0,0,\ldots,0,1_{\sigma},0,\ldots,0)\cdot D_F$. Donc la donnée d'une filtration décroissante exhaustive séparée de $D_F$ par des sous-$F\otimes_{\Q}E$-modules $(\mathrm{Fil}^i D_F)_i$ (pas forcément libres) équivaut à la donnée, pour tout $i\in \mathbb{Z}$ et tout $\sigma\colon F \into E$, d'un sous-$F\otimes_{F,\sigma}E$-espace vectoriel $\mathrm{Fil}^i D_{F,\sigma}$ de $D_{F,\sigma}$ qui vérifie les deux conditions suivantes:
\medbreak
\begin{itemize}
\item[(i)] pour tout $i\in \mathbb{Z}$ et tout $\sigma \in S$ on a $\mathrm{Fil}^{i+1} D_{F,\sigma} \subset \mathrm{Fil}^i D_{F,\sigma}$;
\item[(ii)] pour tout $\sigma \in S$ on a
\[
\bigcup_{i\in \mathbb{Z}} \mathrm{Fil}^i D_{F,\sigma} = D_{F,\sigma} \quad \mbox{et} \quad \bigcap_{i\in \mathbb{Z}} \mathrm{Fil}^i D_{F,\sigma} = 0.
\] 
\end{itemize}

Soit $(\mathrm{Fil}^i D_{F,\sigma})_{i,\sigma}$ une telle filtration. On définit
\[
t_H(D_F) = \sum_{i\in \mathbb{Z}} \sum_{\sigma\colon F\into E} i\, \mathrm{dim}_F(\mathrm{Fil}^i D_{F,\sigma}/\mathrm{Fil}^{i+1} D_{F,\sigma}).
\]
La filtration est dite \textit{admissible} si 
\medbreak
\begin{itemize}
\item[(i)] $t_H(D_F) = t_N(D)$;
\item[(ii)] $t_H(D_F')\leq t_N(D')$ pour tout $F_0\otimes_{\Q}E$-sous-module $D'$ stable par $\varphi$ et muni de la filtration induite.
\end{itemize} 

\begin{conj}\label{congettura}
Les deux conditions suivantes sont équivalentes:
\medbreak
\begin{itemize}
\item[(i)] La représentation $\rho_{\vec{d}}\otimes \pi$ admet une norme $G$-invariante, i.e. une norme $p$-adique telle que $\|gv\| = \|v\|$ pour tout $g\in G$ et $v \in \rho_{\vec{d}}\otimes \pi$.
\item[(ii)] Il existe un $\varphi$-module $D$ de rang $2$ tel que 
\[
\mathrm{WD}(\varphi,D)^{ss} = (r,V)
\] 
et une filtration admissible $(\mathrm{Fil}^i D_{F,\sigma})_{i,\sigma}$, avec $i \in \mathbb{Z}$ et $\sigma \in S$, sur $D_F$ telle que
\[
\mathrm{Fil}^i D_{F,\sigma}/\mathrm{Fil}^{i+1} D_{F,\sigma} \neq 0 \quad  \Leftrightarrow \quad i\in \{-d_{\sigma}-1, 0 \}. 
\]
\end{itemize} 
\end{conj}

L'implication $(i) \Rightarrow (ii)$ de la Conjecture \ref{congettura} a été démontrée dans \cite[Corollary 3.3]{bs}. Plus précisément, on a le résultat suivant.

\begin{prop}\label{propfin}
Considérons les quatre conditions suivantes:
\medbreak
\begin{itemize}
\item[(i)] La représentation $\rho_{\vec{d}}\otimes \pi$ admet une norme $G$-invariante, i.e. une norme $p$-adique telle que $\|gv\| = \|v\|$ pour tout $g\in G$ et $v \in \rho_{\vec{d}}\otimes \pi$.
\item[(ii)] Il existe un $\varphi$-module $D$ de rang $2$ tel que 
\[
\mathrm{WD}(\varphi,D)^{ss} = (r,V)
\] 
et une filtration admissible $(\mathrm{Fil}^i D_{F,\sigma})_{i,\sigma}$, avec $i \in \mathbb{Z}$ et $\sigma \in S$, sur $D_F$ telle que
\[
\mathrm{Fil}^i D_{F,\sigma}/\mathrm{Fil}^{i+1} D_{F,\sigma} \neq 0 \quad  \Leftrightarrow \quad i\in \{-d_{\sigma}-1, 0 \}. 
\]
\item[(iii)] Il existe un $\varphi$-module $D$ de rang $2$ tel que 
\[
(\varphi^f)^{ss} = \begin{bmatrix} {\alpha} & {0} \cr {0} & {\beta} \end{bmatrix}
\] 
et une filtration admissible $(\mathrm{Fil}^i D_{F,\sigma})_{i,\sigma}$, avec $i \in \mathbb{Z}$ et $\sigma \in S$, sur $D_F$ telle que
\[
\mathrm{Fil}^i D_{F,\sigma}/\mathrm{Fil}^{i+1} D_{F,\sigma} \neq 0 \quad  \Leftrightarrow \quad i\in \{-d_{\sigma}-1, 0 \}. 
\]
\item[(iv)] les inégalités suivantes sont vérifiées:
\medbreak
\begin{align}
val_F(\alpha^{-1})+val_F(p\beta^{-1})+\sum_{\sigma \in S}d_{\sigma} = 0;  \label{ineg1}\\
val_F(p\beta^{-1})+\sum_{\sigma \in S}d_{\sigma} \geq 0. \label{ineg2}
\end{align}
\end{itemize} 
Alors on a les implications et équivalences suivantes:
\[
(i) \Rightarrow (ii) \Leftrightarrow (iii) \Leftrightarrow (iv). 
\]
\begin{proof}
L'implication (i) $\Rightarrow$ (iv) découle de \cite[Lemma 7.9]{pas}
\medbreak
L'équivalence (ii) $\Leftrightarrow$ (iii) est immédiate.
\medbreak
L'équivalence (iii) $\Leftrightarrow$ (iv) est une conséquence de \cite[Proposition 3.2]{bs}.
\end{proof}
\end{prop}

\begin{rem}
{\rm
Si $F = \Q$ et $r$ n'est pas scalaire, la Conjecture \ref{congettura} est vraie (\cite{bb}).}
\end{rem}

D'après la Proposition \ref{propfin} il est clair que pour avoir une réponse positive à la conjecture de Breuil et Schneider pour $\GL(F)$ il suffit de montrer l'implication $(iv) \Rightarrow (i)$. Le cas $\alpha \in \Oe^{\times}$ (resp. $\beta \in \Oe^{\times}$) est facile.  Plus précisément on a la proposition suivante. 

\begin{prop}\label{proposizionefacile}
Supposons que \eqref{ineg1} et \eqref{ineg2} soient satisfaites, et supposons  $\alpha \in \Oe^{\times}$ (resp. $\beta \in \Oe^{\times}$). Alors la représentation $\rho_{\vec{d}}\otimes \pi$ admet une norme $G$-invariante. 
\begin{proof}
\begin{itemize}
\item[$\bullet$] Supposons $\alpha \in \Oe^{\times}$ et notons 
\[
\chi = \mathrm{nr}(\alpha^{-1}) \otimes \mathrm{nr}(p\beta^{-1}) \prod_{\sigma \in S} \sigma^{d_{\sigma}}
\]  
On définit $I^{cont}(\chi)$ comme l'induite continue du caractère $\chi$, i.e. l'espace des $\phi\colon G \to E$ continues et telles que
\[
\forall b \in P \ \mbox{et} \ g \in G, \quad \phi(bg) = \chi(b) \phi(g),  
\]
l'action de $G$ étant la translation usuelle à droite sur les fonctions. D'après \eqref{ineg1} et \eqref{ineg2} le caractère $\chi$ est à valeurs entières ce qui implique que la norme définie par
\[ 
\forall \phi  \in I^{cont}(\chi), \quad  \|\phi \| = \sup_{g\in P\backslash G} \phi(g)
\] 
est une norme $G$-invariante sur $I^{cont}(\chi)$. Or on a une injection $G$-équivariante évidente
\[
\rho_{\vec{d}} \otimes \pi \into I^{cont}(\chi),
\]
d'où le résultat. 
\item[$\bullet$] En utilisant l'entrelacement
\[
\mathrm{Ind}_P^G (\mathrm{nr}(\alpha^{-1}) \otimes \mathrm{nr}(p\beta^{-1})) = \mathrm{Ind}_P^G (\mathrm{nr}(\beta^{-1}) \otimes \mathrm{nr}(p \alpha^{-1})),
\]
un argument analogue au précédent permet de conclure.
\end{itemize} 
\end{proof}
\end{prop}

Supposons maintenant $\alpha,\beta \notin \Oe^{\times}$. Sous certaines conditions, le Théorème \ref{principe} permet de donner une réponse positive à la Conjecture \ref{congettura}.

\begin{cor}
Supposons que \eqref{ineg1} et \eqref{ineg2} soient satisfaites, et supposons que  $\vec{d} = (d_{\sigma})_{\sigma \in S}$ satisfait les conditions (i) et (ii) du Théorème \ref{principe}. Alors la représentation $\rho_{\vec{d}}\otimes \pi$ admet une norme $G$-invariante. 
\begin{proof}
Il suffit de montrer que $\rho_{\vec{d}}\otimes \pi$ possède un $\Oe$-réseau séparé et stable sous l'action de $G$. Or, d'après la Proposition \ref{indotte} on dispose d'un isomorphisme $G$-équivariant
\[
\rho_{\vec{d}}\otimes \pi \simeq  \frac{\ind \rho_{\vec{d}}}{(T-a_p)(\ind \rho_{\vec{d}})}
\]
où $a_p = \alpha^f + \beta^f$. Donc, si le vecteur $\vec{d}$ vérifie les conditions (i) et (ii) du Théorème \ref{principe} alors la Proposition \ref{libero} s'applique  et l'on peut en déduire que la représentation $\rho_{\vec{d}}\otimes \pi$ possède un réseau séparé et stable sous l'action de $G$, d'où le résultat.  
\end{proof}
\end{cor}

\end{document}